\documentclass{article}
\usepackage{graphicx} 

\newcommand{\mytitle}{Extending the Veblen Function}
\newcommand{\myauthor}{Jayde Sylvie Massmann, Adrian Wang Kwon}

\font \rfont = cmr12 at 15pt
\font \authorfont = cmr12 at 12.5pt

\title{\rfont \mytitle}
\author{\authorfont \myauthor}
\date{\rfont \today}

\usepackage{hyperref}
\usepackage{comment}
\usepackage{proof}

\usepackage[utf8]{inputenc}
\usepackage{xcolor}
\usepackage{biblatex}
\usepackage{csquotes}
\usepackage{amsfonts}
\usepackage{amsmath}
\usepackage{cleveref}
\usepackage{amssymb}
\usepackage{amsthm}
\usepackage{parskip}
\usepackage[a4paper, total={6in, 8in}]{geometry}
\usepackage{titling}
\usepackage{tikz}
\usetikzlibrary{cd}
\usetikzlibrary{matrix,arrows}
\usetikzlibrary{decorations.pathreplacing,calc,arrows.meta}
\usepackage[normalem]{ulem}
\usepackage[english]{babel}

\usepackage{blindtext}

\theoremstyle{definition}
\newtheorem{definition}{Definition}[section]

\theoremstyle{plain}
\newtheorem{theorem}[definition]{Theorem}

\theoremstyle{plain}
\newtheorem{conjecture}[definition]{Conjecture}

\theoremstyle{plain}
\newtheorem{proposition}[definition]{Proposition}

\theoremstyle{plain}
\newtheorem{corollary}[definition]{Corollary}

\theoremstyle{remark}
\newtheorem{remark}{Remark}

\theoremstyle{remark}
\newtheorem{example}{Example}

\theoremstyle{remark}

\theoremstyle{plain}
\newtheorem{lemma}[definition]{Lemma}

\theoremstyle{plain}
\newtheorem{claim}[definition]{Claim}

\newcommand{\at}[2]{\begin{matrix}#1\\#2\end{matrix}}

\begin{document}

\maketitle

\begin{abstract}
This paper serves to define an extension, which we call \textit{dimensional Veblen}, of Oswald Veblen's system of ordinal functions below the large Veblen ordinal. This is facilitated by iterating derivatives of ordinal functions along multidimensional array structures, and can be viewed as the ``maximal'' natural extension of the Veblen functions. We then construct an ordinal notation based on it, and provide a conversion algorithm from Buchholz's function below the Bachmann-Howard ordinal.
\end{abstract}

\section{Introduction}
Ordinal numbers were introduced by Georg Cantor in 1883, intended as an extension to the natural numbers. The motivation was to generalize enumeration of well-ordered sets to infinitely large collections -- they also proved useful in the proof of the Cantor-Bendixson theorem. Ordinals were later recast by von Neumann with the following simple definition, making them a versatile tool for modern set theory: ordinals were transitive sets well-ordered by $\in$. For example, natural numbers could be encoded as ordinals by representing them recursively as the set of their predecessors. One could then compare ordinals by representing the less-than relation as the elementhood relation, so that $0 < 1$ since $\emptyset \in \{\emptyset\}$. A simple way of representing sufficiently small ordinals was introduced by Cantor:

\begin{definition}
Say an ordinal $\alpha$ has Cantor normal form (CNF) $\omega^{\alpha_0} + \omega^{\alpha_1} + \cdots + \omega^{\alpha_n}$ iff $\alpha = \omega^{\alpha_0} + \omega^{\alpha_1} + \cdots + \omega^{\alpha_n}$, where $\omega$ is the least infinite ordinal, and $\alpha_0 \geq \alpha_1 \geq \cdots \geq \alpha_n$. $\alpha_0$ is called the degree of $\alpha$ (although this terminology is not very widespread).
\end{definition}

\begin{theorem}
Every ordinal has a Cantor normal form.
\end{theorem}

\begin{proof}
Let $\alpha$ be an arbitrary ordinal. Proceed by transfinite induction (more accurately, well-founded ordinal induction) on $\alpha$ -- one can find a proof in \cite{cantor}, but we give one anyways for the sake of completeness. The case $\alpha = 0$ is trivial, with $n = -1$ (i.e. the Cantor normal form for zero is empty). Therefore suppose $\alpha > 0$, and that all ordinals less than $\alpha$ have a Cantor normal form. There exists a least ordinal $\zeta$ so that $\omega^\zeta > \alpha$, and $\zeta$ must be a successor. As such, $\zeta$ has a direct predecessor -- call it $\delta$ -- and $\delta$ is precisely the maximal ordinal so that $\omega^\delta \leq \alpha$. Now there are unique ordinals $\beta, \gamma$ so that $\alpha = \omega^\delta \cdot \beta + \gamma$ and $\gamma < \omega^\delta$. We clearly must have $\beta < \omega$, since $\beta \geq \omega$ would imply $\omega^\zeta \leq \omega^\delta \cdot \beta \leq \alpha$ (contradicting $\alpha < \omega^\zeta$) and so we can use the inductive assumption that $\gamma$ has a Cantor normal form to get a Cantor normal form for $\alpha$.
\end{proof}

The limit of Cantor normal form is $\varepsilon_0$, the small Cantor ordinal. This is the least ordinal $\alpha$ with degree equal to itself, as opposed to degree less than itself. After the Cantor normal form theorem, further ordinal representation systems were introduced, most notably Veblen's function $\varphi$ -- cf. \cite{veblen}. It was first defined by Oswald Veblen in 1908, and the Klammersymbolen (lit. \textit{bracket symbols}) were defined by Kurt Schütte in 1945 as a variadic (possibly transfinitary) extension -- cf. \cite{schutte}. The following is a heavily simplified version of the binary system, where $\mathrm{Ord}$ represents the class of ordinals, and $\mathrm{AP}$ represents the class of nonzero ordinals which can't be reached from below via addition, e.g. $1$ and $\omega$ (formally, $\alpha \in \mathrm{AP}$ iff $\alpha > 0$ and, for all $\beta, \gamma < \alpha$ we have $\beta + \gamma < \alpha$).

\begin{definition}
Define by transfinite recursion functions $\varphi_\alpha: \mathrm{Ord} \to \mathrm{AP}$, where $\varphi_0(\delta) = \omega^\delta$ and, for $\alpha > 0$, $\varphi_\alpha$ enumerates the class $\{\xi: \forall \gamma < \alpha (\varphi_\gamma(\xi) = \xi)\}$. $\varphi$ can be considered as a function $\mathrm{Ord}^2 \rightarrow \mathrm{AP}$, writing $\varphi(\alpha,\beta)$ for $\varphi_\alpha(\beta)$.
\end{definition}

The lemma which permits this type of transfinitely recursive fixed-point taking is Veblen's fixed point lemma, which can be considered as an analogue of Brouwer's fixed point theorem. Say a function $f: \mathrm{Ord} \to \mathrm{Ord}$ is normal if $\alpha < \beta$ implies $f(\alpha) < f(\beta)$, and $f$ is continuous in the order topology on $\mathrm{Ord}$ -- equivalently, $f$ commutes with suprema, i.e. $f(\sup X) = \sup f''X$ whenever $X$ is a set of ordinals. The normal functions are precisely generalizations of, e.g. continuous monotone functions on the reals.

\begin{lemma}
Assume $f$ is normal. Then there is a proper class of fixed points of $f$, and the function enumerating them is also normal.
\end{lemma}

\begin{proof}
To show this, let $\alpha$ be any ordinal, and construct the sequence $\langle \alpha_i: i < \omega \rangle$ by $\alpha_0 = \alpha$, and $\alpha_{i + 1} = f(\alpha_i)$. Let $\beta = \sup\{\alpha_i: i < \omega\}$. Then

\begin{equation}
\begin{split}
f(\beta) & = f(\sup\{\alpha_i: i < \omega\}) \\ & = \sup\{f(\alpha_i): i < \omega\} \\ & = \sup\{\alpha_{i+1}: i < \omega\} \\ & = \beta
\end{split}
\end{equation}

Therefore, $\beta$ is a fixed point of $f$ greater than or equal to $\alpha$, so the class of fixed points is unbounded and therefore proper. Now let $f'$ be the function enumerating the fixed points of $f$. That $f'$ is increasing is obvious, while that it is normal follows from the fact that, for any limit ordinal $\delta$, if $\delta' = \sup\{f'(\gamma): \gamma < \delta\}$:

\begin{equation}
\begin{split}
f(\delta') & = f(\sup\{f'(\gamma): \gamma < \delta\}) \\ & = \sup\{f(f'(\gamma)): \gamma < \delta\} \\ & = \sup\{f'(\gamma): \gamma < \delta\} \\ & = \delta'
\end{split}
\end{equation}

and so $f(\delta) = \delta'$.
\end{proof}

Veblen's function and CNF were catalysts for even stronger ordinal representation systems, notably the projection functions used in ordinal analysis, e.g. in \cite{buchholz3}. These ``collapse'' large cardinals\footnote{Some systems instead use nonrecursive ordinals such as admissible ordinals instead, but countability of the ordinals involved may make some of the function's properties harder to prove.} to large countable ordinals.

In this paper we introduce an extension of the Veblen functions, which allow one to define ordinals up to the famous Bachmann-Howard ordinal without impredicative means (i.e. reference to nonrecursive or uncountable ordinals). The extension contained herein puts arrays on the bottom row, with structures like

\begin{equation}
\varphi\left(\at{1}{(1,0)}\right)
\end{equation}

as the least instance of this. This corresponds to the large Veblen ordinal, or the least $E$-number. Since our extension allows for taking fixed points over lengths of rows, we refer to this extension as ``dimensional Veblen''. More intuition for this naming is provided within our closing remarks.

We prove the equivalence of our system with the original one, where the two domains overlap. Furthermore, we give a conversion algorithm from an initial segment of Buchholz's function, one of the simplest projection functions, to our system. This implies that the limit of our system is at least the Bachmann-Howard ordinal.

For the rest of this paper, we restrict ourselves to working with ordinals below $\omega_2$, the initial ordinal of $\aleph_2$. This is because most of our work deals with large, countable (in fact, recursive) ordinals, but usage of base-$\omega_1$ polynomials will be crucial in conversion between dimensional Veblen and Buchholz's function. The primary reason why we use this technical restriction in the first place is to ensure that we don't have to work with proper classes. Namely, $\mathrm{Ord}$ can, in the following, be taken to mean $\omega_2$.

\section{Definition}
Numerous auxiliary functions are necessary for the upcoming definitions. We present these auxiliary functions, and then use them to define the extension of the Veblen function. Naturally, we also provide intuition. Note that, in the following, when two or more cases of a definition are both satisfied, the first one listed takes priority. 

\begin{definition}
We call a relation $R$ injective iff, for all $x, y, z, w \in \mathrm{fld}(R)$, $(x,y) \in R$ and $(z,w) \in R$ implies either $x = z$ or $y \neq w$. A function $f$ (when considered as its graph $\{(x,y): f(x) = y\}$) is injective in this sense iff it is injective in the usual sense, however, we will primarily be considering injective non-function-like relations, such as $\{(2,2), (2,1)\}$.

For a set $X$, we let the index-derived set of $X$ consist of injective relations on $\mathrm{Ord} \times X$ whose domain is nonempty yet finite, and does not contain zero. The weak index-derived set is defined analogously, except we permit the domain to include zero.

We define inductively $A_n$ for $n < \omega$ by $A_0 = \{\emptyset\}$ and letting $A_{n+1}$ be the index-derived set of $A_n$. Then $A$ is the union of all the $A_n$, and $A' = A \setminus A_0$. The sequence $B_n$ and their union $B$ is defined analogously, replacing index-derivation with weak index-derivation.
\end{definition}

Index-derivation provides a relatively simple method of indexing elements within a possibly multidimensional, nested array. The entire original system of Schütte's Klammersymbolen is obviously isomorphic, in some natural sense, to the set $A_2$, by replacing the bracket symbol

\begin{equation}
\left(\at{a_1 & a_2 & \cdots & a_n}{b_1 & b_2 & \cdots & b_n}\right)
\end{equation}

with the function $\{(a_m, \underline{b_m}): 1 \leq m \leq n\}$. Here, and in the subsequent work, the following abbreviations are used:

\begin{equation}
\begin{split}
\underline{\alpha} & = \{(\alpha,\emptyset)\} \\ (\alpha_i)_{i < n, n < \omega} & = \{(\alpha_i, \underline{i}): i < n < \omega\} \\ \left(\at{\alpha_i}{X_i} \right)_{i < n, n < \omega} & = \{(\alpha_i, X_i): i < n < \omega\}
\end{split}
\end{equation}

\begin{definition}
Assume $X \in A$. Then we let $a(X)$ be the range of $X$, i.e. $\{X': \exists \alpha \in \mathrm{Ord} ((\alpha, X') \in X)\}$. Also, for $X, Y \in A$, let $c_X(Y)$ be the unique $\alpha$ so that $(\alpha, Y) \in X$ if $Y \in a(X)$ (this is the preimage of $Y$ under $X$) if it exists, and else $c_X(Y) = 0$.

We then simultaneously define functions $b: A' \to A$, $\bar{b}: A' \to A$, and a relation $\prec$ on $A$ like so:

For $X \in A$, $b(X)$ is the unique $X' \in a(X)$ so that, for no $X'' \in a(X)$ do we have $X' \prec X''$. Dually, $\bar{b}(X)$ is the unique $X' \in a(X)$ so that, for no $X'' \in a(X)$ do we have $X'' \prec X'$.

\begin{enumerate}
    \item If $Y = \emptyset$, then $X \prec Y$ is false.
    \item Else if $X = \emptyset$, then $X \prec Y$ is true.
    \item Else if $b(X) \neq b(Y)$, then $X \prec Y$ iff $b(X) \prec b(Y)$.
    \item Else if $c_X(b(X)) \neq c_Y(b(Y))$, then $X \prec Y$ iff $c_X(b(X)) < c_Y(b(Y))$.
    \item Else, $X \prec Y$ iff $X \setminus \{(c_X(b(X)), b(X))\} \prec Y \setminus \{(c_Y(b(Y)), b(Y))\}$.
\end{enumerate}
\end{definition}

Intuitively, $\prec$ compares arrays by their first or most significant point of difference, and then $b$ and $\bar{b}$ return the largest and smallest, respectively, subarrays of an array.

We give some basic results, all of which (other than (2), which requires somewhat heavy case classification) are quite elementary.

\begin{lemma}\;
\begin{enumerate}
    \item $c$ is well-defined.
    \item $b$ and $\bar{b}$ are well-defined.
    \item For all $X \in A'$, $b(X) \prec X$.
    \item $\prec$ is a strict linear order.
    \item For all $X \in A'$, either $\bar{b}(X) = b(X)$ or $\bar{b}(X) \prec b(X)$.
\end{enumerate}
\end{lemma}

\begin{proof}
To facilitate this and following proofs, we introduce the rank function $\rho: A \to \mathbb{N}$, by letting $\rho(X)$ be the unique $n$ so that $X \in A_n$.

(1) We want to prove that, for all $X, Y$, if there is some $\alpha$ so that $(\alpha, Y) \in X$, then there is precisely one such $\alpha$. If there weren't, let $\alpha_0 \neq \alpha_1$ and $(\alpha_0, Y), (\alpha_1, Y) \in X$. Since $X$ is injective, either $\alpha_0 = \alpha_1$ or $Y \neq Y$. Neither of these options is possible.

(2) This is very easy to verify, and follows from trichotomy and linearity, which we just proved.

(3) We prove $b(X) \prec X$ for all $X$ so that $X \neq \emptyset$ (equiv. $\rho(X) > 0$). We proceed by induction on $\rho(X)$.

\begin{enumerate}
    \item \textit{Base case.} Assume $\rho(X) = 1$. Then $b(X) \in a(X) = \{\emptyset\}$, therefore $b(X) = \emptyset$. Since $\rho(X) \neq 0$, we have $X \neq \emptyset$ and so $\emptyset \prec X$.
    \item \textit{Induction step.} Assume $\rho(X) = n+1$ for $n > 0$, and $b(Y) \prec Y$ whenever $\rho(Y) \leq n$. Note that then $\rho(b(X)) = n$ and $\rho(b(b(X))) = n-1$, and so the third case must hold. Therefore $b(X) \prec X$ iff $b(b(X)) \prec b(X)$, which holds by the inductive hypothesis.
\end{enumerate}

(4) First, we show that $X \prec X$ never holds. We proceed by induction along $\subseteq$, which is valid because $X$ must be finite.

\begin{enumerate}
    \item \textit{Base case.} Assume $X = \emptyset$. Then $X \prec X$ is obviously false.
    \item \textit{Induction step.} Assume $X \neq \emptyset$, and we have that $Y \prec Y$ is false whenever $Y \subset X$. The first four clauses can't be activated, so we have $X \prec X$ iff $X \setminus \{(c_X(b(X)), b(X))\} \prec X \setminus \{(c_X(b(X)), b(X))\}$. This is impossible by the inductive hypothesis.
\end{enumerate}

We now show asymmetry, namely that $X \prec Y$ implies $Y \prec X$ is false. We do this by induction on $\max(\rho(X), \rho(Y))$.

\begin{enumerate}
    \item \textit{Base case}. Assume $\max(\rho(X), \rho(Y)) = 0$. Then $X = Y = \emptyset$, and so the hypothesis $X \prec Y$ can't happen in the first place.
    \item \textit{Induction step.} Assume $\max(\rho(X), \rho(Y)) = n+1$ for $n > 0$, and $X' \prec Y'$ implies $Y' \prec X'$ is false whenever $\rho(X'), \rho(Y') \leq n$. Assume $b(X) \neq b(Y)$. Then $Y \prec X$ iff $b(Y) \prec b(X)$, which can't happen because $\max(\rho(b(X)), \rho(b(Y))) = n$. Else, assume $c_X(b(X)) \neq c_Y(b(Y))$. Then $Y \prec X$ iff $c_Y(b(Y)) < c_X(b(X))$, which is false because $X \prec Y$ implies $c_X(b(X)) < c_Y(b(Y))$. Else, we proceed by a layered transfinite induction along $\subseteq$ (which, as mentioned previously, is well-founded when restricted to finite sets). Assume that, for all $X' \subset X$ and $Y' \subset Y$, $X' \prec Y'$ implies that $Y' \prec X'$ is false. Then, $Y \prec X$ iff $Y \setminus \{(c_Y(b(Y)), b(Y))\} \prec X \setminus \{(c_X(b(X)), b(X))\}$, which is false. 
\end{enumerate}

We show linearity (more accurately, trichotomy), again by induction on $\max(\rho(X), \rho(Y))$.

\begin{enumerate}
    \item \textit{Base case.} Assume $\max(\rho(X), \rho(Y)) = 0$. Then $X = Y$.
    \item \textit{Induction step.} Assume $\max(\rho(X), \rho(Y)) = n+1$, and, whenever $\rho(X'), \rho(Y') \leq n$, we either have $X' \prec Y'$, $X' = Y'$ or $Y' \prec X'$. We again proceed by case classification to show either $X \prec Y$, $X = Y$ or $Y \prec X$.
    \begin{enumerate}
        \item If $Y = \emptyset$, then we have $Y \prec X$.
        \item If $X = \emptyset$, then we have $X \prec Y$.
        \item If $b(X) \neq b(Y)$, then, since $\max(\rho(b(X)), \rho(b(Y))) = n$, we either have $b(X) \prec b(Y)$ or $b(Y) \prec X$. These are equivalent, respectively, to $X \prec Y$ or $Y \prec X$.
        \item Else if $c_X(b(X)) \neq c_Y(b(Y))$, then we have $X \prec Y$ or $Y \prec X$ since the usual $<$ relation is trichotomous.
        \item Else if $c_X(b(X)) = c_Y(b(Y))$, then perform a transfinite $\subset$-induction argument analogous to the induction step of the proof of asymmetry.
    \end{enumerate}
\end{enumerate}

Transitivity follows from a similar argument, via induction and case classification.

(5) We now verify that, for all $X \neq \emptyset$, either $\bar{b}(X) = b(X)$ or $\bar{b}(X) \prec b(X)$.

\begin{enumerate}
    \item \textit{Base case.} Assume $\rho(X) = 1$. Then $b(X) = \emptyset$ and $\bar{b}(X) = \emptyset$, therefore $b(X) = \bar{b}(X)$.
    \item \textit{Induction step.} Assume $\rho(X) = n+1$ for $n > 0$, and either $\bar{b}(X) = b(X)$ or $\bar{b}(X) \prec b(X)$ whenever $\rho(X) \leq n$. Then, by definition, we have $b(X), \bar{b}(X) \in a(X)$, therefore we can't have $b(X) \prec \bar{b}(X)$. Therefore, the desired result follows from (2). \qedhere
\end{enumerate}
\end{proof}

\begin{lemma}
For ordinals $\alpha, \beta$ we have $\alpha < \beta$ iff $\underline{\alpha} \prec \underline{\beta}$.
\end{lemma}

\begin{proof}
Since $\underline{\alpha}$ is never empty, cases 1 doesn't apply. We have $b(\underline{\alpha}) = b(\underline{\beta}) = \emptyset$, so cases 2 and 3 doesn't apply, but case 4 does apply since $c_{\underline{\alpha}}(b(\underline{\alpha})) = \alpha$ and $c_{\underline{\beta}}(b(\underline{\beta})) = \beta$ and we can, without loss of generality, assume $\alpha \neq \beta$ -- which yields $\underline{\alpha} \prec \underline{\beta}$ iff $\alpha < \beta$, the desired result.
\end{proof}

\begin{definition}
For $X \in B$, we let $d(X) \in A$ be recursively defined as $d(X) = \{(X',d(X'')): (X',X'') \in X \land X' > 0\}$, which recursively erases zeroes from $X$.

For $X \in A$, we let $e(X) \in A$ be defined as $X$ if $c_X(\emptyset) \in \mathrm{Lim} \cup \{0\}$ (i.e. is not a successor), and else as $d(X \setminus \{(c_X(\emptyset),\emptyset)\}\cup\{(\eta,\emptyset)\})$, where $\eta$ is the unique ordinal so that $\eta+1 = c_X(\emptyset)$. In other words, $e(X)$ is obtained from $X$ by subtracting one from the position of $\emptyset$ in $X$. 

For $X \in A$, we recursively define $f(X) \in \mathrm{Ord}$ like so:

\begin{enumerate}
    \item If $X = \emptyset$, $f(X) = 0$.
    \item If $e(X) \neq X$, $f(X) = 1$.
    \item Suppose $c_X(\bar{b}(X))$ is a successor.
    \begin{enumerate}
        \item If $e(\bar{b}(X)) \neq \bar{b}(X)$, then $f(X) = 2$.
        \item Else, $f(X) = f(\bar{b}(X))$.
    \end{enumerate}
    \item Else, $f(X) = c_X(\bar{b}(X))$.
\end{enumerate}

Also, for $X \in A$ and $\alpha \in \mathrm{Ord}$, we define $g(X, \alpha) \in A$ like so:

\begin{enumerate}
    \item If $X = \emptyset$, $g(X, \alpha) = \emptyset$.
    \item If $e(X) \neq X$, $g(X, \alpha) = e(X)$.
    \item Suppose $c_X(\bar{b}(X))$ is a successor. Let $\eta$ be the unique ordinal so that $\eta+1 = c_X(\bar{b}(X))$.
    \begin{enumerate}
        \item If $e(\bar{b}(X)) \neq \bar{b}(X)$, then $g(X, \alpha) = d(X \setminus \{(c_X(\bar{b}(X)), \bar{b}(X))\} \cup \{(\eta,\bar{b}(X)),(\alpha, e(\bar{b}(X)))\})$.
        \item Else, $g(X, \alpha) = d(X \setminus \{(c_X(\bar{b}(X)), \bar{b}(X))\} \cup \{(\eta,\bar{b}(X)),(1,g(\bar{b}(X),\alpha))\})$.
    \end{enumerate}
    \item Else, $g(X, \alpha) = d(X \setminus \{(c_X(\bar{b}(X)), \bar{b}(X))\} \cup \{(\alpha,\bar{b}(X))\})$.
\end{enumerate}

For $X \in A$, define $h(X) \in A$ as $h(X) = \emptyset$ if $X = \emptyset$; $h(X) = X \setminus \{(c_X(\emptyset), \emptyset)\}$ if $\emptyset \in a(X)$; and $h(X) = X$ else.
\end{definition}

$f$ acts as a case classification function, which will tell us whether an array should be evaluated by taking fixed points or limits (2 or $> 2$, respectively). This is reflected in the definition of $\varphi$. Meanwhile, $g$ yields ``fundamental sequences'' for arrays of a sort. And $h$ just removes the last (or first) entry of the array, namely the one indexed by $\emptyset$. An easy observation is the following:

\begin{proposition}
Assume $X \in A'$. Then, for all $\alpha < c_X(\bar{b}(X))$, we have $g(X, \alpha) \prec X$.
\end{proposition}

This is sharp in basically all contexts, although $\alpha$ can be made arbitrarily large if $e(X) \neq X$.

Finally, all the auxiliary functions have been given. We are now able to define $\varphi: A \to \mathrm{Ord}$ like so:

\begin{definition}
Suppose $X \in A$. If $X = \emptyset$, then $\varphi X = 1$. If $a(X) = \{\emptyset\}$, then $\varphi X = \omega^{c_X(\emptyset)}$. Else, $\varphi X = \operatorname{enum}(U)(c_X(\emptyset))$, where:

\begin{enumerate}
    \item If $f(h(X)) > 2$, then $U = \{\alpha: \forall 0 < \beta < f(h(X)) (\alpha = \varphi d(g(h(X), \beta) \cup \{(\alpha, \emptyset)\}))\}$.
    \item If $f(h(X)) = 2$, then $U = \{\alpha: \alpha = \varphi g(h(X), \alpha)\}$.
\end{enumerate}
\end{definition}

\begin{example}
For example, let us take $\varphi\{(1, \{(2, \emptyset)\}), (3, \{(1, \emptyset)\})\}$. To calculate this, one must first calculate $f(h(X))$, which requires one to calculate $h(X)$. Let $X = \{(1, \{(2, \emptyset)\}), (3, \{(1, \\ \emptyset)\})\}$. $a(X) = \{\{(2, \emptyset)\},\{(1, \emptyset)\}\}$, which does not have $\emptyset$ as an element. Therefore, $h(X) = X$. Then, what is $f(h(X)) = f(X)$? Since $X \neq \emptyset$ and $e(X) = X$, cases 1 and 2 do not apply. Therefore, evaluation requires calculation of $c_X(\bar{b}(X))$. We have $\bar{b}(X) = \{(1,\emptyset)\}$ since $\underline{1} \prec \underline{2}$ and so $c_X(\bar{b}(X)) = 3$. This is a successor, so case 3 applies. $e(\bar{b}(X)) = \emptyset \neq \{(1,\emptyset)\}$, so case 3a applies and $f(X) = 2$. Then, one has to calculate $g(X, \alpha)$. Since $c_X(\bar{b}(X))$ is again a successor, and $e(\bar{b}(X)) \neq \bar{b}(X)$, case 3a applies. Therefore, $g(X, \alpha) = d(X \setminus \{(c_X(\bar{b}(X)), \bar{b}(X))\} \cup \{(\eta,\bar{b}(X)),(\alpha, e(\bar{b}(X)))\}) = d(\{(1, \{(2,\emptyset)\}), (2, \{(1, \emptyset)\}), (\alpha, \emptyset)\})$. Therefore, $\varphi\{(1, \{(2, \emptyset)\}), (3, \{(1, \emptyset)\})\}$ is the first fixed point of $\alpha \mapsto \varphi \{(1, \{(2,\emptyset)\}), (2, \{(1, \emptyset)\}), (\alpha, \emptyset)\}$. It turns out that, by an easy inductive argument, that $\varphi \{(1, \{(2,\emptyset)\}), (2, \{(1, \emptyset)\}), (\alpha, \emptyset)\} = \varphi(1, 2, \alpha)$ -- in particular, $\varphi\{(1, \{(2, \emptyset)\}), (3, \{(1, \emptyset)\})\} = \varphi(1, 3, 0)$.
\end{example}

\begin{example}
Let us take $\varphi({}^2_\omega)$ = $\varphi\{(2,\{(\omega,\emptyset)\})\}$. $f(h(X)) = f(X)$ again, and case 3b applies, as $e(\{(\omega,\emptyset)\}) = \{(\omega,\emptyset)\}$, so $f(X) = \omega$. Then, calculate $g(X,\alpha)$. Since $e(\overline{b}(X)) = \overline{b}(X)$, case 3b applies. Therefore, $g(X,\alpha) = d(\{(1,\underline{\omega}),(1,\underline{\alpha})\})$. And since $f(X) = \omega > 2$, $\varphi\{(2,\{(\omega,\emptyset)\})\}$ is the first common fixed point of $\alpha\mapsto\{(1,\{(\omega,\emptyset)\}),(1,\{(\gamma,\emptyset)\}),(\alpha,0)\} = \varphi({}^1_\omega\;{}^1_\gamma\;{}^\alpha_0)$ for $\gamma<\omega$, or the limit of $\varphi({}^1_\omega\;{}^1_\alpha)$ as $\alpha$ tends to $\omega$.
\end{example}

\section{An ordinal notation}

We provide a ``constructive'' ordinal notation system (a way of representing ordinals in a computable way), associated to the dimensional Veblen function. We then gauge the strength of the dimensional Veblen function, and its associated ordinal notation.

The usage of corner brackets in the following is to distinguish between the ordinals and functions on them, and terms representing them, although we'll drop them later. Of course, it is possible to encode these terms as natural numbers, but we will stick with the approach utilizing formal strings.

\begin{definition}
We let $\mathcal{A}, \mathcal{T}$ be sets of formal strings defined recursively like so:

\begin{enumerate}
    \item $\ulcorner 0 \urcorner \in \mathcal{T}$.
    \item If $a, b \in \mathcal{T}$, then $\ulcorner a + b \urcorner \in \mathcal{T}$.
    \item If $a \in \mathcal{A}$, then $\ulcorner \varphi a \urcorner \in \mathcal{T}$.
    \item If $n < \omega$, $\alpha_1, \alpha_2, \cdots, \alpha_n \in \mathcal{T}$ and $a_1, a_2, \cdots, a_n \in \mathcal{A}$, then $(\alpha_1 @ a_1, \alpha_2 @ a_2, \cdots, \alpha_n @ a_n) \in \mathcal{A}$.
    \item $\epsilon$ is an abbreviation for $()$.
\end{enumerate}

We determine the subterms of $t \in \mathcal{T} \cup \mathcal{A}$, denoted $\mathbf{S}(t)$, like so:

\begin{enumerate}
    \item $\mathbf{S}(a) = \{a\}$ if $a \in \mathcal{T}$.
    \item $\mathbf{S}((\alpha_0 @ a_0, \alpha_1 @ a_1, \cdots, \alpha_n @ a_n)) = \{\alpha_0, \alpha_1, \cdots, \alpha_n\} \cup \bigcup_{i \leq n} \mathbf{S}(a_i)$
\end{enumerate}

In particular, $\mathbf{S}(\epsilon) = \emptyset$.

We now simultaneously define relations $a \prec b$ on $\mathcal{A}$ and $s < t$ on $\mathcal{T}$. $\prec$ is defined as (the smallest relation satisfying) the following:

\begin{enumerate}
    \item If $a = \epsilon$, then $a \prec b$ iff $a \neq b$.
    \item If $b = \epsilon$, then $a \prec b$ is false.
    \item Else, let $a = (\alpha_1 @ a_1, \alpha_2 @ a_2, \cdots, \alpha_n @ a_n)$ and $b = (\beta_1 @ b_1, \beta_2 @ b_2, \cdots, \beta_m @ b_m)$ for $n, m > 0$. Recursively let $a'$ and $b'$ be the $\prec$-maximal elements of $\{a_1, a_2, \cdots, a_n\}$ and $\{b_1, b_2, \cdots, b_n\}$, respectively. Also, let $\iota_a$ and $\iota_b$ be the unique natural numbers so that $a' = a_{\iota_a}$ and $b' = b_{\iota_b}$. Then:
    \begin{enumerate}
        \item If $a' \neq b'$, then $a \prec b$ iff $a' \prec b'$.
        \item If $\alpha_{\iota_a} \neq \beta_{\iota_b}$, then $a \prec b$ iff $\alpha_{\iota_a} \prec \beta_{\iota_b}$.
        \item Else, $a \preceq b$ iff $(\alpha_i@a_i)_{i \neq \iota_a \land i < n} \prec (\beta_i@b_i)_{i \neq \iota_b \land i < m}$.
    \end{enumerate}
\end{enumerate}

Obviously, this is an implementation of $\prec$ (on $A$). And for $<$:

\begin{enumerate}
    \item If $a = 0$, $a < b$ iff $a \neq b$.
    \item If $b = 0$, $a < b$ is false.
    \item Suppose $a = a' + a''$ for $a', a'' \in \mathcal{T}$.
    \begin{enumerate}
        \item Suppose $b = b' + b''$ for $b', b'' \in \mathcal{T}$.
        \begin{enumerate}
            \item If $a' \neq b'$, $a < b$ iff $a' < b'$.
            \item If $a' = b'$, $a < b$ iff $a'' < b''$.
        \end{enumerate}
        \item Suppose $b = \varphi B$ for $B \in \mathcal{A}$.
        \begin{enumerate}
            \item If $a' \neq b$, $a < b$ iff $a' < b$.
            \item If $a' = b$, $a < b$ is false.
        \end{enumerate}
    \end{enumerate}
    \item Suppose $a = \varphi A$ for some $A$.
    \begin{enumerate}
        \item Suppose $b = b' + b''$ for $b', b'' \in \mathcal{T}$.
        \begin{enumerate}
            \item If $a \neq b'$, $a < b$ iff $a < b'$
            \item If $a = b'$, $a < b$ is true.
        \end{enumerate}
        \item If $b = \varphi B$ for some $B$, $a < b$ iff either of the following hold:
        \begin{enumerate}
            \item For all $k \in \mathbf{S}(B)$, we have $a < k$.
            \item $A \prec B$ and, for all $k \in \mathbf{S}(A)$, we have $k < b$.
        \end{enumerate}
    \end{enumerate}
\end{enumerate}
\end{definition}

Obviously, $<$ is a comparison algorithm on (terms representing) ordinals. The last clause is analogous to the fact that, unlike some other ordinal notations, comparison $\varphi(\alpha_1, \beta_1) < \varphi(\alpha_2, \beta_2)$ for principal terms in Veblen normal form is non-lexicographical, and instead has three cases:

\begin{enumerate}
    \item $\alpha_1 < \alpha_2$ and $\beta_1 < \varphi(\alpha_2, \beta_2)$.
    \item $\alpha_1 = \alpha_2$ and $\beta_1 < \beta_2$.
    \item $\alpha_2 < \alpha_1$ and $\varphi(\alpha_2, \beta_2) < \beta_1$.
\end{enumerate}

but our case is obviously slightly more complex as all positions and subarrays in the arrays must also be considered.

However, $<$ gives incorrect results such as $\ulcorner 0 \urcorner < \ulcorner 0 + 0 \urcorner$. Therefore, we make our final definition -- the notions of principality and standardness:

\begin{definition}
We let $\mathcal{PT}$ be the set of principal terms, where $s$ is principal iff it is of the form $\varphi a$ for some $a \in \mathcal{A}$.
\end{definition}

Next, we must define standardness; which yields a \textit{unique} representation for every ordinal below the limit. We don't necessarily have unique representations already -- for example, letting $\$1$ abbreviate $\varphi \epsilon$, we want that $\varphi(\$1 @ (\$1 @ (\$1 @ \epsilon)))$ (which corresponds to $E_0$) is standard, but $\varphi(\$1 @ \varphi(\$1 @ (\$1 @ \\ (\$1 @ \epsilon))))$ (which would correspond to $\varphi(1 @ E_0) = E_0$) isn't -- this is because the former is a less pathological representation of the same ordinal, and allowing duplicates could allow for a very erratic, ill-founded structure. Checking standardness of sums involves ensuring they're written in Cantor normal form, and checking standardness of principal terms involves hereditarily making sure all the subterms are standard, and that they aren't fixed points of the function they've been put into. First, we must define a few more auxilary functions on $\mathcal{A}$.

\begin{definition}
Assume $X \in \mathcal{A}$. Let $X = (\alpha_1 @ a_1, \alpha_2 @ a_2, \cdots, \alpha_n @ a_n)$. Then, let $p(X)$ be the unique $k$ so that, for no $i \neq k$ do we have $a_i \prec a_k$. This plays an analogue to $\bar{b}$. Also we define $m(X) \in \mathcal{T}$ as $\ulcorner 0 \urcorner$ if $n = 0$ (i.e. $X = \epsilon$), and $\alpha_{p(X)}$ else.

Assume $a \in \mathcal{A}$. Then $q(X, a) = (\alpha_i@a_i)_{i < n \land a \prec a_i}$.
\end{definition}

Now we can finally define standardness:

\begin{definition}
We define a set $\mathcal{OT}$ like so:

\begin{enumerate}
    \item $0 \in \mathcal{OT}$.
    \item $\alpha+\beta \in \mathcal{OT}$ iff all of the following hold:
    \begin{enumerate}
        \item $\alpha \in \mathcal{PT} \cap \mathcal{OT}$.
        \item $\beta \in \mathcal{OT} \setminus \{0\}$.
        \item $\beta = \alpha$ or $\beta < \alpha$.
    \end{enumerate}
    \item $\varphi X \in \mathcal{OT}$ iff $\mathbf{S}(X) \subset \mathcal{OT}$ and one of the following holds:
    \begin{enumerate}
        \item $X = \epsilon$.
        \item $m(X) \notin \mathcal{PT}$.
        \item For some $a \in \mathbf{S}(q(X,p(X)))$, we have $m(X) \leq a$.
        \item $q(X,p(X)) \preceq q(X',0)$, where $X'$ is the unique element of $\mathcal{A}$ such that $\varphi X' = m(X)$.
    \end{enumerate}
\end{enumerate}
\end{definition}

The (standard) principal terms are intended to correspond to the additively principal ordinals (elements of $\mathrm{AP}$) since all other terms can be represented from zero and addition. This section culminates in a proof that the ``limit'' of dimensional Veblen, equal to the supremum of $\varepsilon_0$, $E_0$ (the large Veblen ordinal), $\varphi(1@(1@(1,0)))$, and so on, is at least the Bachmann-Howard ordinal $\eta_0$ -- and we conjecture the two are equal. This would imply that the elements of $\mathcal{PT} \cap \mathcal{OT}$ correspond precisely to $\mathrm{AP} \cap \eta_0$. We do this by connecting our dimensional Veblen function to Buchholz's $\psi$-functions, defined in \cite{buchholz}.

In the definitions that follow, $\mathcal{J} = \varepsilon_{\Omega+1}$, where $\Omega$ is the (initial ordinal of) the least uncountable cardinal $\aleph_1$. It's known that Buchholz's $\psi$-functions have the following relationships with the usual Veblen function and Klammersymbolen (where $\alpha, \beta, \gamma, \cdots$ are sufficiently small):

\begin{enumerate}
    \item $\psi_0(\alpha) = \varphi(0, \alpha)$.
    \item $\psi_0(\Omega^\alpha (1+\beta)) = \varphi(\alpha, \beta)$, for $\beta > 0$.
    \item $\psi_0(\Omega^{\Omega^n \alpha_n + \cdots + \Omega \alpha_1 + \alpha_0} (1+\beta)) = \varphi(\alpha_n, \cdots, \alpha_1, \alpha_0, \beta)$, for $\beta, n > 0$.
\end{enumerate}

In particular:

\begin{tabular}{|cccccl|}
     \hline
     \textbf{Ordinal} && \textbf{in $\psi$} && \textbf{in $\psi$} &  \\
     \hline
     $\varepsilon_0$ &=& $\psi_0(\Omega)$ &=& $\psi_0(\psi_1(0))$ & (small Cantor ordinal) \\
     $\zeta_0$ &=& $\psi_0(\Omega^2)$ &=& $\psi_0(\psi_1(\psi_1(0)))$ & (Cantor ordinal) \\
     $\Gamma_0$ &=& $\psi_0(\Omega^\Omega)$ &=& $\psi_0(\psi_1(\psi_1(\psi_1(0))))$ & (Feferman-Schütte ordinal) \\
     $\varphi(1,0,0,0)$ &=& $\psi_0(\Omega^{\Omega^2})$ &=& $\psi_0(\psi_1(\psi_1(\psi_1(0)2)))$ & (Ackermann ordinal) \\
     $\varphi\left(\at{1}{\omega}\right)$ &=& $\psi_0(\Omega^{\Omega^\omega})$ &=& $\psi_0(\psi_1(\psi_1(\psi_1(1))))$ & (small Veblen ordinal) \\
     $E_0$ &=& $\psi_0(\Omega^{\Omega^\Omega})$ &=& $\psi_0(\psi_1(\psi_1(\psi_1(\psi_1(0)))))$ & (large Veblen ordinal) \\
     \hline
\end{tabular}

We provide a conversion algorithm between Buchholz's $\psi$ functions, sub-$\mathcal{J}$, to the collection of standard $\varphi$-terms. This algorithm also requires numerous auxiliary functions.

When we decompose an ordinal into the form $\xi+\Omega^\gamma \delta$, it is always assumed that $\xi$ is of the form $\Omega^{\gamma+1} \eta$ and that $0 < \delta < \Omega$. Note that $\eta$ may be zero.

\begin{definition}
For $\alpha < \mathcal{J}$, we define $s(\alpha)$ like so:

\begin{enumerate}
    \item If $\alpha = 0$, then $s(\alpha) = \emptyset$.
    \item If $\alpha > 0$, let $\alpha = \xi + \Omega^\beta \gamma$. Then, $s(\alpha) = \{\beta, \gamma\} \cup s(\xi)$.
\end{enumerate}

\begin{remark}
$s(\alpha)$ is the set of all exponents and coefficients in $\alpha$ when written in base-$\Omega$ CNF. Note that this is not hereditary -- in the latter case we do not include $s(\beta)$ or $s(\gamma)$.
\end{remark}

For $\alpha, \beta < \mathcal{J}$, we recursively define our trichotomic case-classification function $k(\alpha,\beta)$ like so:

\begin{enumerate}
    \item Suppose $\alpha < \Omega$.
    \begin{enumerate}
        \item If $\alpha < \beta$, then $k(\alpha,\beta) = -1$.
        \item If $\alpha = \beta$, then $k(\alpha,\beta) = 0$.
        \item If $\alpha > \beta$, then $k(\alpha,\beta) = 1$.
    \end{enumerate}
    \item Suppose $\alpha \geq \Omega$. Then, write $\alpha$ as $\xi + \Omega^\gamma \delta$.
    \begin{enumerate}
        \item Suppose that, for all $\rho \in s(\alpha)$, we have $k(\rho, \beta) = -1$. Then $k(\alpha,\beta) = -1$.
        \item Suppose that, for all $\rho \in s(\xi)$, we have $k(\rho, \beta) = -1$; and either:
        \begin{enumerate}
            \item $\gamma = \beta$ and $\delta = 1$.
            \item $k(\gamma, \beta) = -1$ and $\delta = \beta$.
        \end{enumerate}
        Then $k(\alpha, \beta) = 0$.
        \item Otherwise, $k(\alpha,\beta) = 1$.
    \end{enumerate}
\end{enumerate}
\end{definition}

This next function distinguishes between three cases for the conversion of larger ordinals into smaller ordinals, which can then be used when converting $\psi_0$ into $\varphi$. For arbitrary ordinals $\rho_1 \geq \rho_2$, let $\rho_1 -- \rho_2$ be the unique ordinal so that $\rho_2 + (\rho_1 -- \rho_2) = \rho_1$.

\begin{definition}
For $\alpha < \mathcal{J}$, we define $t(\alpha)$ like so.

\begin{enumerate}
    \item If $\alpha = 0$, $t(\alpha) = 0$.
    \item Suppose $\alpha > 0$. Then, let $\alpha = \xi + \Omega^\beta \gamma$. Set $\lambda = \psi_0(\xi) -- 1$, and $u = k(\beta, \lambda)$.
    \begin{enumerate}
        \item If $u = -1$, $\rho = \lambda$. 
        \item If $u = 0$, $\rho = 1$.
        \item If $u = 1$, $\rho = 0$.
    \end{enumerate}
    Finally, let $t(\alpha) = \Omega \beta + (\rho + \gamma -- 1)$.
\end{enumerate}
\end{definition}

\begin{lemma}
$t$ is injective.
\end{lemma}

\begin{proof}
Assume $t(\alpha) = t(\beta)$. It's easy to verify that $t(\alpha) = 0$ iff $\alpha = 0$ (not just if) so we get both $\alpha = 0$ and $\beta = 0$ if just either holds. Else, write $\alpha = \xi_1 + \Omega^{\beta_1} \gamma_1$ and $\beta = \xi_2 + \Omega^{\beta_2} \gamma_2$. Let $\lambda_1, \lambda_2, u_1, u_2, \rho_1, \rho_2$ be defined as expected. Since have $\rho_1 + \gamma_1 -- 1, \rho_2 + \gamma_2 -- 1 < \Omega$ we get $\beta_1 = \beta_2$ and $\rho_1 + \gamma_1 = \rho_2 + \gamma_2$. It's relatively easy to verify via case classification that $\xi_1 = \xi_2$. From this it follows that $\gamma_1 = \gamma_2$, by combining $\rho_1 = \rho_2$, a consequence of $\xi_1 = \xi_2$, with $\rho_1 + \gamma_1 = \rho_2 + \gamma_2$, and so $\alpha = \beta$.
\end{proof}

\begin{definition}
For $\alpha < \mathcal{J}$, we define $V(\alpha)$ like so.
\begin{enumerate}
    \item If $\alpha = 0$, then $V(\alpha) = \emptyset$.
    \item If $\alpha > 0$, then let $\alpha = \xi+\Omega^\beta \gamma$. Then, $V(\alpha) = V(\xi) \cup \{(\gamma, V(\beta))\}$.
\end{enumerate}
\end{definition}

In the following, for simplicity, we introduce some notation. Firstly, we use the traditional way of writing summation (in, say number theory): $\sum_{i = 1}^n \alpha_i$ is meant to be interpreted as $\alpha_1 + \alpha_2 + \cdots + \alpha_n$. This is not $\alpha_n + \alpha_{n-1} + \cdots + \alpha_1$: this sum would instead be denoted $\sum_{i = n}^1 \alpha_i$, although we have no usage for ``inverse sums'' of such a form. Note that this allows us to succinctly define $V$ as$\rho_1 + \gamma_1 = \rho_2 + \gamma_2$

\begin{equation}
V\left(\sum_{i = 1}^n \Omega^{\beta_i} \gamma_i\right) = \{(\gamma_i, V(\beta_i)): 1 \leq i \leq n\}
\end{equation}

when $0 < \gamma_i < \Omega$ for all $1 \leq i \leq n$, and $\beta_1 > \beta_2 > \cdots > \beta_n$. This can be used to obtain an analogue of Lemma 2.4:

\begin{proposition}
For ordinals $\alpha, \beta < \mathcal{J}$ we have $\alpha < \beta$ iff $V(\alpha) \prec V(\beta)$.
\end{proposition}

\begin{proof}
If $\alpha = 0$, then $V(\alpha) = \emptyset$, and so $V(\alpha) \prec V(\beta)$ iff $V(\beta) \neq \emptyset$ iff $\beta > 0$. If $\beta = 0$, then $V(\beta) = \emptyset$, and so $V(\alpha) \prec V(\beta)$ is impossible: similarly, $\alpha < \beta$ is impossible. Now we can assume $\alpha, \beta$ are nonzero. Let $\alpha = \xi_1 + \Omega^{\gamma_1} \delta_1$ and $\beta = \xi_2 + \Omega^{\gamma_2} \delta_2$. Then $V(\alpha) = V(\xi_1) \cup \{(\delta_1, V(\gamma_1))\}$ and $V(\alpha) = V(\xi_2) \cup \{(\delta_2, V(\gamma_2))\}$. We proceed by well-founded induction along the lexicographical order on ordinals, so assume that $\alpha' < \beta'$ iff $V(\alpha') \prec V(\beta')$ whenever either $\alpha' < \alpha$ or $\alpha' = \alpha$ and $\beta' < \beta$. If $\xi_1 = \xi_2 = 0$, we get $b(V(\alpha)) = V(\gamma_1)$ and $b(V(\beta)) = V(\gamma_2)$, so $\gamma_1 \neq \gamma_2$ implies $V(\alpha) \prec V(\beta)$ iff $V(\gamma_1) \prec V(\gamma_2)$, which is equivalent to $\gamma_1 < \gamma_2$. And this is itself equivalent to $\alpha < \beta$ due to $\xi_1 = \xi_2 = 0$ and $0 < \delta_1, \delta_2 < \Omega$. If $\gamma_1 = \gamma_2$ then $V(\alpha) \prec V(\beta)$ iff $c_{V(\alpha)}(b(V(\alpha))) < c_{V(\beta)}(b(V(\beta)))$ (since we can assume, without loss of generality, that $\alpha \neq \beta$), and $c_{V(\alpha)}(b(V(\alpha))) = \delta_1$. Similarly, $c_{V(\beta)}(b(V(\beta))) = \delta_2$, so $V(\alpha) \prec V(\beta)$ iff $\delta_1 < \delta_2$ iff $\alpha < \beta$. The rest of the argument (when $\xi_1$ or $\xi_2$ is nonzero) is completely analogous -- for example, if $\xi_1 = 0$ and $\xi_2 > 0$ then  $b(V(\beta)) = b(V(\xi_2))$, so $V(\alpha) \prec V(\beta)$ iff $V(\gamma_1) \prec b(V(\xi_2))$ -- one can then fully break down $\xi_2$ further into base-$\Omega$ Cantor normal form, and then proceed using the inductive hypothesis.
\end{proof}

\begin{theorem}
Let $\mathcal{D} = \{\alpha: \alpha \in C_0(\alpha)\} \cap \mathcal{J}$. Then, for all $\alpha \in \mathcal{D}$, $\psi_0(\alpha) = \varphi V(t(\alpha))$. Here $C$ and $\psi$ are again from \cite{buchholz}.
\end{theorem}

Note that the upcoming proof implicitly uses some results from \cite{buchholz} without mention of them.

\begin{proof}
We proceed by transfinite induction.

\begin{enumerate}
    \item Assume $\alpha = 0$. Then $t(\alpha) = 0$ and so $V(t(\alpha)) = V(0) = \emptyset$. We clearly have $\varphi \emptyset = 1$, and $\psi_0(0) = 1$ as well.
    \item Assume $\alpha = \alpha'+1$, for some $\alpha'$. Clearly, we have $\alpha' \in \mathcal{D}$ too, which will be useful later since it'll permit us to apply the inductive hypothesis. If $\alpha' = 0$ then $\lambda = 0$ so $u = k(0, 0) = 0$ and $\rho = 1$. Therefore $t(\alpha) = 1$ and $V(t(\alpha)) = \{(1,\emptyset)\}$. Then $\varphi V(t(\alpha)) = \omega^1 = \omega$ and $\psi_0(\alpha) = \psi_0(1) = \omega$. Else, let $\alpha' = \xi+\Omega^\gamma \delta$ where $\xi = \Omega^{\gamma+1} \eta$. 
    \begin{enumerate}
        \item If $\gamma = 0$, then $\alpha = \xi + \Omega^\gamma (\delta+1)$. Then $t(\alpha) = \rho+\delta+1-1$. If $\alpha' < \Omega$, then $\alpha < \Omega$ too and so $\xi = \lambda = u = 0$ and $\rho = 1$. Therefore, $t(\alpha) = \delta+1$ and $\varphi V(t(\alpha)) = \varphi (V(0) \cup \{(\delta+1,V(0))\}) = \varphi \{(\delta+1,\emptyset)\} = \omega^{\delta+1} = \omega^\alpha$. By $\alpha \in C_0(\alpha)$, we have $\psi_0(\alpha) = \omega^\alpha$.
        
        Else, $\xi > 0$. We have $\lambda = \psi_0(\xi)$. Since $\gamma = 0$, we have $\gamma < \lambda < \Omega$ and so $u = -1$, therefore $\rho = \lambda$. Then $t(\alpha)  = \lambda+\delta+1$. Therefore, $V(t(\alpha)) = V(\lambda+\delta+1) = \{(\lambda+\delta+1, \emptyset)\}$ and so $\varphi V(t(\alpha)) = \omega^{\lambda+\delta+1}$. Note that $t(\alpha') = \lambda + \delta -- 1$. It follows that $\varphi V(t(\alpha)) = \varphi V(t(\alpha')) \omega$. Then, since $\alpha' \in \mathcal{D}$, we get $\psi_0(\alpha) = \psi_0(\alpha') \omega$, and so the desired result follows from the inductive hypothesis.
        
        \item Else, $\gamma > 0$. Then $\alpha = \Omega (\Omega^{\gamma+1-1} \eta + \Omega^{\gamma-1} \delta)+1$, and so $\lambda = \psi_0(\alpha')$. This yields $u = k(0,\lambda) = -1$. Therefore, $t(\alpha) = \lambda+1$ and so $V(t(\alpha)) = \{(\lambda+1, \emptyset)\}$, giving $\varphi V(t(\alpha)) = \omega^{\lambda+1} = \omega^{\psi_0(\alpha')+1} = \psi_0(\alpha') \omega$. Since $\alpha' \in C_0(\alpha')$, we have $\psi_0(\alpha) = \psi_0(\alpha'+1) = \psi_0(\alpha') \omega = \varphi V(t(\alpha))$.
    \end{enumerate}
    \item Assume $\alpha$ is a limit ordinal. Let $\alpha = \psi_{i_1}(\mu_1) + \psi_{i_2}(\mu_2) + \cdots + \psi_{i_k}(\mu_k)$, where $\psi_{i_1}(\mu_1) \geq \psi_{i_2}(\mu_2) \geq \cdots \geq \psi_{i_k}(\mu_k)$, $\{i_1, i_2, \cdots, i_k\} \subseteq \{0,1\}$, $\mu_j \in C_{i_j}(\mu_j)$ for $1 \leq j \leq k$, and either $i_k = 1$ or $\mu_k > 0$. We split this up into two cases, whether or not $\alpha$ is further additively principal (equivalently, whether or not $k = 1$):
    \begin{enumerate}
        \item $\alpha$ is not additively principal. If $i_1 = 0$ then $t(\alpha) = \alpha$ and so $\varphi V(t(\alpha)) = \varphi \{(\alpha,\emptyset)\} = \omega^\alpha$, and $\psi_0(\alpha) = \omega^\alpha$ too since $\mathcal{D} \cap \Omega = \varepsilon_0$.

        Else, let $\ell$ be such that, for all $1 \leq j \leq k$, we have $i_j = 1$ iff $1 \leq j \leq \ell$. This is well-defined since $i_j = 0$ for all $j$ implies $\alpha < \Omega$, and we can't have $i_j = 0$ yet $i_{j'} = 1$ for $j < j'$. Set

        \begin{equation}
        \begin{split}
        \alpha_1 & = \sum_{j = 1}^\ell \psi_{i_j}(\mu_j) \\
        \alpha_2 & = \sum_{j = \ell+1}^k \psi_{i_j}(\mu_j) \\
        \end{split}
        \end{equation}

        so that $\alpha = \alpha_1 + \alpha_2$, $\alpha_1 \geq \Omega$ and $\alpha_2 < \Omega$.

        \begin{enumerate}
            \item If $\alpha_2 = 0$ (i.e. $\ell = k$), then $\alpha$ is a multiple of $\Omega$. Let $\tau$ (for ``term'') be $k -- \min\{j: \psi_{i_{j+1}}(\mu_{j+1}) = \psi_{i_k}(\mu_k)\}$, the amount of repetitions of the last term. Also let $\delta, \xi$ be the degree and coefficient, respectively, of the last term. This yields

            \begin{equation}
            \alpha = \sum_{j = 1}^{k-\tau} \psi_{i_j}(\mu_j) + \Omega^\delta \xi \tau
            \end{equation}
            
            Also note that $\xi$ is additively principal. Then $\lambda = \psi_0\left(\sum_{j = 1}^{k-\tau} \psi_{i_j}(\mu_j)\right)$ and $u = k(\delta, \lambda)$. We have a few cases to consider:
        
            \begin{enumerate}
                \item If $\delta < \lambda < \Omega$ then $u = -1$ so $\rho = \lambda$ and $t(\alpha) = \Omega \delta + \lambda + \xi \tau$. Then $\varphi V(t(\alpha)) = \varphi (\{(\delta, V(1)), (\lambda + \xi \tau, \emptyset)\})$. 
        
                \begin{claim}
                For all $\eta, \delta$, we have $\varphi \{(\eta, V(1)), (\delta, \emptyset)\} = \varphi(\eta, \delta)$.
                \end{claim}
            
                \begin{proof}
                We use well-founded ordinal induction on $\eta$. Naturally, we need to evaluate $f(h(\{(\eta, \{(1,\emptyset)\}), (\delta, \emptyset)\}))$. This is precisely $f(\{(\eta, \{(1,\emptyset)\})\})$. If $\eta$ is a limit ordinal, then this is equal to $\eta > 2$, so
            
                \begin{equation}
                \begin{split}
                U & = \{\xi: \forall \beta < \eta (\xi=\varphi d(g(h(\{(\eta, \{(1,\emptyset)\}), (\delta, \emptyset)\}),\beta) \cup \{(\xi, \emptyset)\}))\} \\ & = \{\xi: \forall \beta < \eta (\xi=\varphi d(g(\{(\eta, \{(1,\emptyset)\})\},\beta) \cup \{(\xi, \emptyset)\}))\} \\ & = \{\xi: \forall 0 < \beta < \eta (\xi=\varphi \{(\beta, \{(1,\emptyset)\}), (\xi, \emptyset)\})\}
                \end{split}
                \end{equation}
            
                By the inductive hypothesis, $\varphi \{(\beta, \{(1,\emptyset)\}), (\xi, \emptyset)\} = \varphi(\beta,\xi)$ so $U = \{\xi: \forall \beta < \eta (\xi=\varphi(\beta,\xi))\}$ and then $\varphi \{(\eta, \{(1,\emptyset)\}), (\delta, \emptyset)\} = \operatorname{enum}(U)(\delta) = \varphi(\eta, \delta)$.
            
                Meanwhile, if $\eta$ is successor, then $f(\{(\eta, \{(1,\emptyset)\})\}) = 2$ so
            
                \begin{equation}
                \begin{split}
                U & = \{\xi: \xi = \varphi g(\{(\eta, \{(1,\emptyset)\})\}, \xi)\} \\ & = \{\xi: \xi = \varphi \{(\dot{\eta},\{(1,\emptyset)\}),(\xi, \emptyset)\}\}
                \end{split}
                \end{equation}
            
                where $\dot{\eta}$ is the unique ordinal so that $\dot{\eta}+1 = \eta$. Note that $\dot{\eta} = \eta'$ iff $\eta < \omega$.
            
                By the inductive hypothesis, $\varphi \{(\dot{\eta},\{(1,\emptyset)\}),(\xi, \emptyset)\} = \varphi(\dot{\eta}, \xi)$ so $U = \{\xi: \xi = \varphi(\dot{\eta}, \xi)\}$ and $\varphi \{(\eta, \{(1,\emptyset)\}), (\delta, \emptyset)\} = \operatorname{enum}(U)(\delta) = \varphi(\dot{\eta}+1, \delta) = \varphi(\eta, \delta)$.
                \end{proof}
            
                And so $\varphi V(t(\alpha)) = \varphi(\delta, \lambda + \xi \tau) = \sup\{\varphi(\delta, \lambda + \xi (\tau -- 1) + \theta): \theta < \xi\}$. Meanwhile
        
                \begin{equation}
                \begin{split}
                \psi_0(\alpha) & = \psi_0\left(\sum_{j = 1}^{k-\tau} \psi_{i_j}(\mu_j) + \Omega^\delta \xi \tau\right) \\ & = \sup\left\{\psi_0\left(\sum_{j = 1}^{k-\tau} \psi_{i_j}(\mu_j) + \Omega^\delta \xi (\tau-1) + \Omega^\delta \theta \right): \theta < \xi\right\} \\ & = \sup\left\{\varphi V\left(t\left(\sum_{j = 1}^{k-\tau} \psi_{i_j}(\mu_j) + \Omega^\delta \xi (\tau-1) + \Omega^\delta \theta\right)\right): \theta < \xi\right\} \\ & = \sup\{\varphi V(\Omega \delta + \lambda + \xi (\tau-1) + \theta): \theta < \xi\} \\ & = \sup\{\varphi (\{(\delta, V(1)), (\lambda + \xi (\tau-1) + \theta, \emptyset)\}): \theta < \xi\} \\ & = \sup\{\varphi(\delta, \lambda + \xi (\tau-1) + \theta): \theta < \xi\}
                \end{split}
                \end{equation}
                
                \item If $\lambda \leq \delta < \Omega$ then $u = 0, 1$ so $\rho = 1, 0$ and $t(\alpha) = \Omega \delta + \xi \tau$. So $\varphi V(t(\alpha)) = \varphi(\delta, \xi \tau)$, and a similar argument as to the previous case applies for $\psi_0(\alpha)$. 
                
                \item If $\delta \geq \Omega$ and, for all $\rho \in s(\delta)$ we have $k(\rho, \lambda) = -1$, then $u = -1$. Namely, let $\delta = \Omega^{\beta+1} \eta + \Omega^\beta \gamma$. Then $t(\alpha) = \Omega^{1+\beta+1} \eta + \Omega^{1+\beta} \gamma + \lambda + \xi \tau$ and $V(t(\alpha)) = V(\Omega^{1+\beta+1} \eta) \cup \{(\gamma, V(1+\beta)), (\lambda + \xi \tau, \emptyset)\}$. Meanwhile,
        
                \begin{equation}
                \begin{split}
                \psi_0(\alpha) & = \psi_0\left(\sum_{j = 1}^{k-\tau} \psi_{i_j}(\mu_j) + \Omega^{\Omega^{\beta+1} \eta + \Omega^\beta \gamma} \xi \tau\right) \\ & = \operatorname{enum}\left\{\varrho: \forall \theta < \gamma \left(\varrho = \psi_0\left(\sum_{j = 1}^{k-\tau} \psi_{i_j}(\mu_j) + \Omega^{\Omega^{\beta+1} \eta + \Omega^\beta \theta} \varrho\right)\right)\right\}(\xi \tau) \\ & = \operatorname{enum}\left\{\varrho: \forall \theta < \gamma \left(\varrho = \varphi V\left(t\left(\sum_{j = 1}^{k-\tau} \psi_{i_j}(\mu_j) + \Omega^{\Omega^{\beta+1} \eta + \Omega^\beta \theta} \varrho\right)\right)\right)\right\} \\ & \; (\xi \tau) \\ & = \operatorname{enum}\{\varrho: \forall \theta < \gamma (\varrho = \varphi V(\Omega^{1+\beta+1} \eta + \Omega^{1+\beta} \theta + \lambda + \varrho))\}(\xi \tau) \\ & = \operatorname{enum}\{\varrho: \forall \theta < \gamma (\varrho = \varphi (V(\Omega^{1+\beta+1} \eta) \cup \{(\theta, V(1+\beta)), (\lambda + \varrho, \emptyset)\})\} \\ & \; (\xi \tau)
                \end{split}
                \end{equation}
        
                Which is the same as $\varphi (V(\Omega^{1+\beta+1} \eta) \cup \{(\gamma, V(1+\beta)), (\lambda + \xi \tau, \emptyset)\})$.
        
                \item Else, $u = 0, 1$ so $\rho = 1, 0$. Then $t(\alpha) = \Omega \delta + \xi \tau$ and $V(t(\alpha)) = \{(\beta, V(1)), (\xi \tau, \emptyset)\}$. Therefore, (like before) we have $\varphi V(t(\alpha)) = \varphi(\beta, \xi \tau)$. A similar argument to the previous case applies to finding $\psi_0(\alpha)$.
            \end{enumerate}
        \item If $\alpha_2 > 0$, we have $\lambda = \psi_0(\alpha_1)$ and $u = k(0, \lambda) = -1$ so $\rho = \lambda$. Therefore $t(\alpha) = \psi_0(\alpha_1) + \alpha_2$. This means $V(t(\alpha)) = \{(\psi_0(\alpha_1) + \alpha_2,\emptyset)\}$ and $\varphi V(t(\alpha)) = \omega^{\psi_0(\alpha_1) + \alpha_2}$. Meanwhile, $\psi_0(\alpha) = \psi_0(\alpha_1 + \alpha_2) = \psi_0(\alpha_1) \psi_0(\alpha_2) = \psi_0(\alpha_1) \omega^{\alpha_2} = \omega^{\psi_0(\alpha_1) + \alpha_2}$.
        \end{enumerate}
            
        \item $\alpha$ is additively principal. Let $\mu = \mu_1$ and $i = i_1$. We will, naturally, distinguish between $i = 0$ and $i = 1$.
        \begin{enumerate}
            \item Suppose $i = 0$. Then $\alpha < \varepsilon_0$, and so $\xi = \lambda = u = 0$ and $\rho = 1$, yielding $\varphi V(t(\alpha)) = \varphi V(\alpha) = \varphi \{(\alpha, \emptyset)\} = \omega^\alpha = \psi_0(\alpha)$.
            \item Suppose $i = 1$. Since $\mu < \mathcal{J}$, we necessarily have $\mu \in C_1(\mu)$. Let $\mu = \xi + \Omega^\gamma \delta$ and $\xi = \Omega^{\gamma+1} \eta$.
            \begin{enumerate}
                \item If $\gamma = 0$, then $\alpha = \psi_1(\Omega \eta + \delta)$ for $0 < \delta < \Omega$, and $\alpha \in \mathcal{D}$ implies this is equal to $\omega^{\Omega (1 + \eta) + \delta} = \Omega^{1 + \eta} \omega^\delta$, so $t(\alpha) = \Omega (1 + \eta) + \rho + \omega^\delta -- 1$. Since $\xi = 0$, we have $\lambda = 0$ so $u = k(1+\eta, 0) = 1$, yielding $\rho = 0$ and $t(\alpha) = \Omega (1 + \eta) + \omega^\delta$. Then $V(t(\alpha)) = V(\Omega (1 + \eta)) \cup \{(\omega^\delta, \emptyset)\} = \{(1+\eta, \{(1,\emptyset)\}), (\omega^\delta, \emptyset)\}$. As such, we get $\varphi V(t(\alpha)) = \varphi(1+\eta, \omega^\delta)$. Meanwhile $\alpha \in C_0(\alpha)$ and $\alpha = \omega^{\Omega (1 + \eta) + \delta} = \Omega^{1+\eta} \omega^\delta$ give $\psi_0(\alpha) = \varphi(1+\eta, \omega^\delta)$, as desired.

                \begin{equation}
                \begin{split}
                \alpha & = \omega^{\Omega+\mu} \\ & = \Omega \omega^{\Omega^{\gamma+1} \eta} \omega^{\Omega^\gamma \delta} \\ & = \Omega \Omega^{\Omega^{\gamma+1-1} \eta} \Omega^{\Omega^{\gamma-1} \delta} \\ & = \Omega^{1 + \Omega^{\gamma+1-1} \eta + \Omega^{\gamma-1} \delta}
                \end{split}
                \end{equation}

                and so

                \begin{equation}
                \begin{split}
                t(\alpha) & = \Omega (1 + \Omega^{\gamma+1-1} \eta + \Omega^{\gamma-1} \delta) + (\rho+1 -- 1) \\ & = \Omega + \Omega^{\gamma+1} \eta + \Omega^\gamma \delta + (\rho+1 -- 1) \\ & = \Omega + \mu + (\rho+1 -- 1)
                \end{split}
                \end{equation}

                We get $u = k(1 + \Omega^{\gamma+1-1} \eta + \Omega^{\gamma-1} \delta, 0) = 1$ so $\rho = 0$ and $t(\alpha) = \Omega + \mu = \Omega + \Omega^{\gamma+1} \eta + \Omega^\gamma \delta$. If $\eta = 0$ and $\gamma = 1$, then this is $\Omega (1 + \delta)$. Then $V(t(\alpha)) = \{(1+\delta, V(1))\}$ and so $\varphi V(t(\alpha)) = \varphi(1+\delta, 0)$. Meanwhile $\psi_0(\alpha) = \psi_0(\Omega^{1 + \delta})$, which is also equal to $\varphi(1+\delta, 0)$ by $\alpha \in C_0(\alpha)$.
                
                Else if $\eta > 0$ or $\gamma > 1$, then $t(\alpha)$ is just $\Omega^{\gamma+1} \eta + \Omega^\gamma \delta$. This gives $V(t(\alpha)) = V(\Omega^{\gamma+1} \eta) \cup \{(\delta, V(\gamma))\}$. Then

                \begin{equation}
                \begin{split}
                f(h(V(t(\alpha)))) & = f(V(t(\alpha))) \\ & = f(V(\Omega^{\gamma+1} \eta) \cup \{(\delta, V(\gamma))\})
                \end{split}
                \end{equation}

                If $\delta$ is a limit ordinal, then this is $\delta > 2$, so

                \begin{equation}
                \begin{split}
                \varphi V(t(\alpha)) & = \min\{\xi: \forall \beta < f(h(V(t(\alpha)))) (\xi=\varphi d(g(h(V(t(\alpha))),\beta) \cup \{(\xi, \emptyset)\}))\} \\ & = \min\{\xi: \forall \beta < \delta (\xi=\varphi d(d(V(\Omega^{\gamma+1} \eta) \cup \{(\beta,V(\gamma))\}) \cup \{(\xi, \emptyset)\}))\} \\ & = \min\{\xi: \forall \beta < \delta (\xi=\varphi (V(\Omega^{\gamma+1} \eta) \cup \{(\beta,V(\gamma)), (\xi, \emptyset)\}))\} \\ & = \sup\{\min\{\xi: \xi=\varphi (V(\Omega^{\gamma+1} \eta) \cup \{(\beta,V(\gamma)), (\xi, \emptyset)\})\}: \beta < \delta\}
                \end{split}
                \end{equation}

                while $\alpha \in C_0(\alpha)$ coupled with the inductive hypothesis yields

                \begin{equation}
                \begin{split}
                \psi_0(\alpha) & = \psi_0(\Omega^{1 + \Omega^{\gamma+1-1} \eta + \Omega^{\gamma-1} \delta}) \\ & = \sup\{\psi_0(\Omega^{1 + \Omega^{\gamma+1-1} \eta + \Omega^{\gamma-1} \beta}): \beta < \delta\} \\ & = \sup\{\varphi V(t(\Omega^{1 + \Omega^{\gamma+1-1} \eta + \Omega^{\gamma-1} \beta})): \beta < \delta\} \\ & = \sup\{\varphi V(\Omega + \Omega^{\gamma+1} \eta + \Omega^\gamma \beta): \beta < \delta\} \\ & = \sup\{\varphi (V(\Omega^{\gamma+1} \eta) \cup \{(\beta, V(\gamma))\}): \beta < \delta\}
                \end{split}
                \end{equation}

                and a simple analysis of the behaviour of $\varphi$ shows that, for all $\beta < \delta$, we have that $\min\{\xi: \xi=\varphi (V(\Omega^{\gamma+1} \eta) \cup \{(\beta,V(\gamma)), (\xi, \emptyset)\})\}$ is in-between $\varphi (V(\Omega^{\gamma+1} \eta) \cup \{(\beta, V(\gamma))\})$ and $\varphi (V(\Omega^{\gamma+1} \eta) \cup \{(\beta+1, V(\gamma))\})$ -- therefore, the two sets have the same suprema and $\varphi V(t(\alpha)) = \psi_0(\alpha)$.

                Meanwhile, if $\delta$ is successor, there are two further cases: whether or not $\gamma$ is successor.

                If $\gamma$ is a limit then $f(h(V(t(\alpha)))) = f(V(\gamma))$. As such, we have another two cases to consider: if $f(V(\gamma)) = 2$ or $f(V(\gamma)) > 2$. In the first case we get

                \begin{equation}
                \begin{split}
                \varphi V(t(\alpha)) & = \min\{\xi: \xi = \varphi g(h(V(t(\alpha))), \xi)\} \\ & = \min\{\xi: \xi = \varphi g(V(\Omega^{\gamma+1} \eta) \cup \{(\delta, V(\gamma))\}, \xi)\} \\ & = \min\{\xi: \xi = \varphi d(V(\Omega^{\gamma+1} \eta) \cup \{(\dot{\delta},V(\gamma)),(1,g(V(\gamma),\xi))\})\}
                \end{split}
                \end{equation}

                where $\dot{\delta}$ is the unique ordinal so that $\dot{\delta}+1 = \delta$.

                \begin{claim}
                For all $\alpha < \mathcal{J}$, we have $f(V(\alpha)) = 2$ iff $\operatorname{cof}(\alpha) = \Omega$.
                \end{claim}

                \begin{proof}
                In the forwards direction, let $\alpha = \xi + \Omega^\gamma \delta$, where $0 < \delta < \Omega$ -- we claim $\delta$ is a successor, from which the desired result follows. Towards contradiction, if $\delta$ is a limit then $f(V(\alpha)) = f(V(\xi) \cup \{(\delta, V(\gamma))\})$. We have $e(V(\alpha)) = V(\alpha)$ due to $\delta$ being a limit, $\bar{b}(V(\alpha)) = V(\gamma)$ and $c_{V(\alpha)}(\bar{b}(V(\alpha))) = \delta$, yielding $f(V(\alpha)) = \delta > 2$, contradicting $f(V(\alpha)) = 2$. The converse follows by a similar argument.
                \end{proof}

                \begin{remark}
                As such, we're basically classifying whether or not $\gamma$ has countable cofinality. If it does, then $\psi_0(\alpha)$ is just a supremum of smaller inputs, due to $\psi_\mu$ being $\omega$-continuous, and similarly $\varphi V(t(\alpha))$ is just a simultaneous fixed point, which reduces to a supremum of smaller inputs (e.g. the fact that $\varphi(\omega, 0)$ is just the supremum of $\varphi(n, 0)$ for $n < \omega$) by Veblen's fixed point lemma combined with the fact that the intersection of set-many clubs in $\mathrm{Ord}$ is still club.
                \end{remark}

                In this case then, $\operatorname{cof}(\gamma) = \Omega$. Let $t(\gamma) = \Omega + \Omega^{\beta+1} \eta' + \Omega^\beta \delta'$. Necessarily, $\delta'$ is successor since $0 < \delta' < \Omega$ and $\delta'$ being a limit would imply $\operatorname{cof}(\delta') = \omega$ and so $\operatorname{cof}(\gamma) = \omega$, contradicting $\operatorname{cof}(\gamma) = \Omega$. The pattern of the subsequent case classification is quite similar to what has come before in this proof -- we leave the exact details to the reader.
                
                In the latter case (where $\gamma$ has countable cofinality) we get

                \begin{equation}
                \begin{split}
                \varphi V(t(\alpha)) & = \min\{\xi: \forall \beta < f(V(\gamma)) (\xi = \varphi d(g(V(t(\alpha)),\beta) \cup \{(\xi, \emptyset)\})))\} \\ & = \min\{\xi: \forall \beta < f(V(\gamma)) (\xi = \varphi d(V(\Omega^{\gamma+1} \eta) \cup \{(\dot{\delta}, V(\gamma)), \\ & \; (1,g(V(\gamma),\beta)), (\xi, \emptyset)\})))\} \\ & = \sup\{\min\{\xi: \xi = \varphi d(V(\Omega^{\gamma+1} \eta) \cup \{(\dot{\delta}, V(\gamma)), (1,g(V(\gamma),\beta)), \\ & \; (\xi, \emptyset)\})))\}: \beta < f(V(\gamma))\} \\ & = \sup\{\varphi d(V(\Omega^{\gamma+1} \eta) \cup \{(\dot{\delta}, V(\gamma)), (1,g(V(\gamma),\beta)), (1,V(1))\}): \beta \\ & \; < f(V(\gamma))\}
                \end{split}
                \end{equation}

                where $\dot{\delta}$ is the unique ordinal so that $\dot{\delta}+1 = \delta$. Meanwhile

                \begin{equation}
                \begin{split}
                \psi_0(\alpha) & = \psi_0(\Omega^{1 + \Omega^{\gamma+1} \eta + \Omega^\gamma \delta}) \\ & = \psi_0(\Omega^{1 + \Omega^{\gamma+1} \eta + \Omega^\gamma \dot{\delta} + \Omega^\gamma}) \\ & = \psi_0(\Omega^{1 + \Omega^{\gamma+1} \eta + \Omega^\gamma \dot{\delta}} \Omega^{\Omega^\gamma}) \\ & = \sup\{\psi_0(\Omega^{1 + \Omega^{\gamma+1} \eta + \Omega^\gamma \dot{\delta}} \Omega^{\Omega^\beta}): \beta < \gamma\} \\ & = \sup\{\varphi V(t(\Omega^{1 + \Omega^{\gamma+1} \eta + \Omega^\gamma \dot{\delta}} \Omega^{\Omega^\beta})): \beta < \gamma\} \\ & = \sup\{\varphi V(\Omega^{\gamma+1} \eta + \Omega^\gamma \dot{\delta} + \Omega^\beta): \beta < \gamma\} \\ & = \sup\{\varphi d(V(\Omega^{\gamma+1} \eta) \cup \{(\dot{\delta}, V(\gamma)), (1,V(\beta))\}): \beta < \gamma\}
                \end{split}
                \end{equation}

                and the two are obviously equal.

                Else, if $\gamma$ is successor then $f(h(V(t(\alpha)))) = 2$ and so

                \begin{equation}
                \begin{split}
                \varphi V(t(\alpha)) & = \min\{\xi: \xi = \varphi g(h(V(t(\alpha))), \xi)\} \\ & = \min\{\xi: \xi = \varphi (V(\Omega^{\gamma+1} \eta) \cup \{(\dot{\delta}, V(\gamma)),(\xi, V(\dot{\gamma}))\})\}
                \end{split}
                \end{equation}

                where $\dot{\gamma}, \dot{\delta}$ are the unique ordinals so that $\dot{\gamma}+1 = \gamma$ and $\dot{\delta}+1 = \delta$.

                Also,

                \begin{equation}
                \begin{split}
                \psi_0(\alpha) & = \psi_0(\Omega^{1 + \Omega^{\gamma+1-1} \eta + \Omega^{\gamma-1} \delta}) \\ & = \psi_0(\Omega^{1 + \Omega^{\gamma+1-1} \eta + \Omega^{\gamma-1} \dot{\delta}} \Omega^{\Omega^{\gamma-1}})
                \end{split}
                \end{equation}

                We now do one last distinction -- whether or not $\gamma$ is finite. If $\gamma$ is finite then $\dot{\gamma} = \gamma-1$ is itself also successor so, letting $\ddot\gamma$ be so that $\ddot\gamma+1 = \dot{\gamma}$:

                \begin{equation}
                \begin{split}
                \psi_0(\Omega^{1 + \Omega^{\gamma+1-1} \eta + \Omega^{\gamma-1} \dot{\delta}} \Omega^{\Omega^{\gamma-1}}) & = \psi_0(\Omega^{1 + \Omega^\gamma \eta + \Omega^{\dot{\gamma}} \dot{\delta}} \Omega^{\Omega^{\ddot\gamma+1}}) \\ & = \min\{\xi: \xi = \psi_0(\Omega^{1 + \Omega^\gamma \eta + \Omega^{\dot{\gamma}} \dot{\delta}} \Omega^{\Omega^{\ddot\gamma} \xi})\} \\ & = \min\{\xi: \xi = \varphi V(\Omega^{\gamma+1} \eta + \Omega^\gamma \dot{\delta} + \Omega^{\dot{\gamma}} \xi)\}
                \end{split}
                \end{equation}

                This becomes $\min\{\xi: \xi = \varphi (V(\Omega^{\gamma+1} \eta) \cup \{(\dot{\delta}, V(\gamma)), (\xi,V(\dot{\gamma}))\})\}$, which is obviously the same as $\varphi V(t(\alpha))$. A similar argument applies when $\gamma$ is infinite, although note that in this case we have $\gamma -- 1 = \gamma$ and so no have no need for $\ddot\gamma$.

                This concludes the proof. \qedhere
            \end{enumerate}
        \end{enumerate}
    \end{enumerate}
\end{enumerate}
\end{proof}

For example, referring back to the previous table:

\begin{enumerate}
    \item $\varepsilon_0 = \varphi V(t(\Omega)) = \varphi V(\Omega) = \varphi \{(1,\{(1,\emptyset)\})\}$.
    \item $\zeta_0 = \varphi V(t(\Omega^2)) = \varphi V(\Omega 2) = \varphi \{(2,\{(1,\emptyset)\})\}$.
    \item $\Gamma_0 = \varphi V(t(\Omega^\Omega)) = \varphi V(\Omega^2) = \varphi \{(1,\{(2,\emptyset)\})\}$.
    \item $\varphi(1, 0, 0, 0) = \varphi V(t(\Omega^{\Omega^2})) = \varphi V(\Omega^3) = \varphi \{(1,\{(3,\emptyset)\})\}$.
    \item $\mathrm{SVO} = \varphi V(t(\Omega^{\Omega^\omega})) = \varphi V(\Omega^\omega) = \varphi \{(1,\{(\omega,\emptyset)\})\}$.
    \item $\mathrm{LVO} = \varphi V(t(\Omega^{\Omega^\Omega})) = \varphi V(\Omega^\Omega) = \varphi \{(1,\{(1,\{(1,\emptyset)\})\})\}$
\end{enumerate}

The array $\{(\alpha,\{(\beta,\emptyset)\})\}$ can be imagined as a $\alpha$ at position $\beta$, which would be written as

\begin{equation}
\left(\at{\alpha}{\beta}\right)
\end{equation}

in Schütte's original system (which in fact agrees with our use of such matrices). This also has some parallels with Weiermann's $\vartheta$, introduced in \cite{rathjen}, although an explicit description of a conversion algorithm would be significantly more arduous. A list of some equivalences is given below. 
\begin{enumerate}
    \item For $\alpha < \zeta_0$, $\vartheta(\alpha) = \varepsilon_\alpha$,
    \item For $\beta < \Omega, \alpha < \zeta_{\beta+1}$, $\vartheta(\zeta_\beta+(\alpha-1)) = \varepsilon_{\zeta_\beta+\alpha}$. 
    And more generally,
    \item For $\eta < \Gamma_0, \alpha < \varphi(1+\eta+1,0)$, $\vartheta(\Omega \eta + \alpha) = \varphi(1+\eta,\alpha)$,
    \item For $\eta < \Gamma_0, \beta < \Omega, \alpha < \varphi(\eta+1,\beta+1)$, $\vartheta(\Omega \eta +\varphi(1+\eta+1,\beta)+(\alpha-1)) = \varphi(1+\eta,\varphi(1+\eta+1,\beta)+\alpha)$.
    \item $\vartheta(\Omega^2) = \Gamma_0$.
    \item $\vartheta(\Omega \Gamma_0) = \varphi(\Gamma_0,1)$.
\end{enumerate}

In fact, $\vartheta$ has no fixed points and has a nicer correspondence with Veblen's $\varphi$. This is in contrast with the behaviour of Buchholz's $\psi$, which is monotonically non-decreasing but also not injective due to getting ``stuck'', e.g. on the interval $[\varepsilon_0, \Omega]$.

As we mentioned earlier, the array $\{(\alpha,\{(\beta,\emptyset)\})\}$ can be imagined as a $\alpha$ at position $\beta$, and corresponds to the ordinal $\vartheta(\Omega^\beta \alpha)$ . In general, terms below the large Veblen ordinal -- reached at $\vartheta(\Omega^\Omega) = \varphi\left(\at{1}{(1,0)}\right)$ -- correspond to 1-dimensional arrays, terms below $\vartheta(\Omega^{\Omega^2}) = \varphi\left(\at{1}{(1,0,0)}\right)$ correspond to 2-dimensional arrays, and so on, even up to $\omega$ dimensions and beyond. The limit of typical intuition is reached by $\vartheta(\Omega^{\Omega^\Omega}) = \varphi\left(\at{1}{\left(\at{1}{(1,0)}\right)}\right)$, a ``hyperdimension'', or row of dimensions. After this, planes and $\alpha$-dimensional structures of dimensions can be reached, and on to a hyperdimensional structure of dimensions, whence it can be nested infinitely.

Now, we introduce a much more straightforward algorithm of converting from dimensional Veblen below the large Veblen ordinal to Veblen's original functions.

First, we define a modification of the $\varphi$-functions from \cite{veblen}. In the following, subscripts denote positions, not variable modifiers.

\begin{definition}
Define the $\varphi^*$-functions as follows:
\begin{enumerate}
    \item Let $\varphi^*(\alpha_0) = \omega^\alpha$.
    \item Let $\varphi^*(\alpha_0,0_1, \cdots,\alpha_\beta, \cdots,\alpha_\gamma)$ be the $\alpha_0$th simultaneous solution of
    
    \begin{equation}
    \eta = \varphi^*(0_0,0_1, \cdots,\eta_\delta,0_{\delta+1}, \cdots,\alpha_\beta ', \cdots,\alpha_\gamma)
    \end{equation}
    
    for all $\delta < \beta, \alpha_\beta' < \alpha_\beta$.
\end{enumerate}
\end{definition}

Now, we will give the seven classes, in this modernized version of the Veblen function. These are crucial for comparison between Veblen's function and our modification. In the following section, instead of $\varphi^*$, we shall write $\varphi$. Also, from now on we shall omit all the entries $0_\alpha$. Lastly, we will let $\alpha@\beta$ stand for $\at{\alpha}{\beta}$ if we intend for the equations to fit inline.

\begin{enumerate}
    \addtocounter{enumi}{-1}
    \item $\varphi(0) = 1$.
    \item \textit{Class A.} $\varphi(1_\alpha)$ for $\alpha \in \mathrm{Succ}$.
    \item \textit{Class B.} $\varphi(1_\alpha)$ for $\alpha \in \mathrm{Lim}$.
    \item \textit{Class C.} $\varphi(\alpha_0, \cdots)$ for $\alpha_0 > 0$\footnote{This has been changed slightly.}.
    \item \textit{Class D.} $\varphi(\alpha_\beta, \cdots)$ for $\beta \in \mathrm{Succ}$, $\alpha \in \mathrm{Succ}$.
    \item \textit{Class E.} $\varphi(\alpha_\beta, \cdots)$ for $\alpha_\beta \in \mathrm{Lim}$.\footnote{The distinction between classes $E$ and $G$ does not need to be made here.}
    \item \textit{Class F.} $\varphi(\alpha_\beta, \cdots)$ for $\beta \in \mathrm{Lim}$, $\alpha \in \mathrm{Succ}$.
\end{enumerate}

Below $\dot{\alpha}$ denotes the unique ordinal so that $\dot{\alpha}+1 = \alpha$, when $\alpha$ is a successor ordinal, generalizing the notation used in the proof of Theorem 3.10. Now, we can define the fundamental sequences for the modernized system. To calculate $\alpha[n]$, which is the $(1+n)^\mathrm{th}$ element of the sequence for $\alpha$:

\begin{enumerate}
    \item If $\alpha$ is simply $\varphi(\beta_0)$ for $\beta$ a successor, $\alpha[n] = \varphi(\dot{\beta}_0) \cdot n$. If $\beta$ is a limit, $\alpha[n] = \varphi(\beta[n]_0)$.
    \item If $\alpha$ is of class \textit{A}, then let $\alpha = \varphi(1_\beta)$. Then $\alpha[0] = \varphi(1_{\dot{\beta}})$ and $\alpha[n+1] = \varphi(\alpha[n]_{\dot{\beta}})$.
    \item If $\alpha$ is of class \textit{B}, then let $\alpha = \varphi(1_\beta)$. Then $\alpha[n] = \varphi(1_{\beta[n]})$.
    \item If $\alpha$ is of class \textit{C}, then let $\alpha = \varphi(\beta_0, \cdots)$. If $\beta \in \mathrm{Lim}$, $\alpha[n] = \varphi(\beta[n]_0, \cdots)$. Otherwise, let $\alpha' = \varphi(0_0, \cdots)$, and let $\rho = \varphi(\dot{\beta}_0, \cdots)$. Then,
    \begin{enumerate}
        \item If $\alpha'$ is of class \textit{A}, then let $\alpha' = \varphi(1_\gamma)$. Then $\alpha[0] = \rho+1$ and $\alpha[n+1] = \varphi(\alpha[n]_{\dot{\gamma}}, \cdots)$.
        \item If $\alpha'$ is of class \textit{B}, then let $\alpha' = \varphi(1_\gamma)$. Then $\alpha[n] = \varphi((\rho+1)_0,1_{\gamma[n]}, \cdots)$\footnote{In the original paper Veblen writes $\varphi((\rho+1)_{\gamma[n]})$, but we have changed this to be consistent with the case for class \textit{F}.}.
        \item If $\alpha'$ is of class \textit{D}, then let $\alpha' = \varphi(\gamma_\delta, \cdots)$. Then $\alpha[0] = \rho+1$, and $\alpha[n+1] = \varphi(\alpha[n]_{\dot{\delta}}, \cdots)$.
        \item If $\alpha'$ is of class \textit{E}, then let $\alpha' = \varphi(\gamma_\delta, \cdots)$. Then $\alpha[n] = \varphi((\rho+1)_0, \alpha[n]_\delta, \cdots)$.
        \item If $\alpha'$ is of class \textit{F}, then let $\alpha' = \varphi(\gamma_\delta, \cdots)$. Then $\alpha[n] = \varphi((\rho+1)_0,1_{\delta[n]}, \dot{\gamma}_\delta, \cdots)$.
    \end{enumerate}
    \item If $\alpha$ is of class \textit{D}, then let $\alpha = \varphi(\beta_\gamma, \cdots)$. Then $\alpha[0] = 0$, and $\alpha[n+1] = \varphi(\alpha[n]_{\dot{\gamma}}, \dot{\beta}_\gamma, \cdots)$.
    \item If $\alpha$ is of class \textit{E}, then let $\alpha = \varphi(\beta_\gamma, \cdots)$. Then $\alpha[n] = \varphi(\beta[n]_\gamma, \cdots)$.
    \item If $\alpha$ is of class \textit{F}, then let $\alpha = \varphi(\beta_\gamma, \cdots)$. Then $\alpha[n] = \varphi(1_{\gamma[n]}, \dot{\beta}_\gamma, \cdots)$.
\end{enumerate}

If one sets $0[n] = 1[n] = 0$ then one can use Cantor normal form to fill in the fundamental sequences for non-additively principal ordinals, as is done in the Wainer hierarchy.

\begin{lemma}
For any normal function $f$, let $g = \operatorname{enum}\{\alpha: \alpha = f(\alpha)\}$ be its derivative. Veblen's fixed point lemma guarantees that this function exists and is total. Then, for all $\alpha$, we have $g(\alpha+1) = \sup\{f^n(g(\alpha)+1): n < \omega\}$
\end{lemma}

\begin{proof}
We claim first that $\sup\{f^n(g(\alpha)+1): n < \omega\}$ is a fixed point of $f$:

\begin{equation}
\begin{split}
f(\sup\{f^n(g(\alpha)+1): n < \omega\}) & = \sup\{f(f^n(g(\alpha)+1)): n < \omega\} \\ & = \sup\{f^{n+1}(g(\alpha)+1): n < \omega\} \\ & = \sup\{f^n(g(\alpha)+1): n < \omega\}
\end{split}
\end{equation}

And as such, $g(\alpha+1) \leq \sup\{f^n(g(\alpha)+1): n < \omega\}$. One can easily show that $f^n(g(\alpha)+1) \leq g(\alpha+1)$ for all $n$, by induction on $n$. This is because $g$ is strictly increasing (so the base case $n = 0$ is trivial), and that any fixed point of a normal function is also a closure point. From this it follows $\sup\{f^n(g(\alpha)+1): n < \omega\} \leq g(\alpha+1)$ and so $g(\alpha+1) = \sup\{f^n(g(\alpha)+1): n < \omega\}$
\end{proof}

\begin{theorem}
If $n < \omega$, for all $\alpha_0, \alpha_1, \cdots, \alpha_n, \beta_0, \beta_1, \cdots, \beta_n$:

\begin{equation}
\varphi\left(\at{\alpha_0 & \alpha_1 & \cdots & \alpha_n}{\underline{\beta_0} & \underline{\beta_1} & \cdots & \underline{\beta_n}}\right) = \varphi^*((\alpha_0)_{\beta_0}, (\alpha_1)_{\beta_1}, \cdots, (\alpha_n)_{\beta_n})
\end{equation}
\end{theorem}

\begin{proof}
Let $\Phi = \varphi^*((\alpha_0)_{\beta_0},(\alpha_1)_{\beta_1}, \cdots,(\alpha_n)_{\beta_n})$, $X = \left(\at{\alpha_0 & \alpha_1 & \cdots & \alpha_n}{\underline{\beta_0} & \underline{\beta_1} & \cdots & \underline{\beta_n}}\right)$, where $\beta_0 > \beta_1 > \cdots > \beta_n$, and $\Phi' = \varphi X$. We use well-founded ordinal induction on $\Phi$ to prove $\Phi = \Phi'$. Assume that we have proven this for every $\alpha < \Phi$. Then, proceed by case classification:
\begin{enumerate}
    \item If $a(X) = \{\emptyset\}$, i.e. $n = \beta_0 = 0$, then $\varphi(\alpha@0) = \varphi^*(\alpha_0) = \omega^\alpha$.
    \item If $\Phi$ is of class $A$, then $\Phi' = \varphi(1@\beta)$, and is the least fixed point of $\xi \mapsto \varphi(\xi@\dot{\beta})$. And $\Phi[n] = (\lambda \xi.\varphi^*(\xi_{\dot{\beta}}))^n (\varphi^*(1_{\dot{\beta}}))$, and so $\Phi$ is the least fixed point of $\xi \mapsto \varphi^*(\xi_{\dot{\beta}})$. These are equal by the inductive hypothesis.
    \item If $\Phi$ is of class $B$, then $\Phi' = \varphi(1@\beta)$ for a limit ordinal $\beta$, and is the least simultaneous fixed point of $\xi \mapsto \varphi(1@\gamma, \xi@0)$ for $\gamma < \beta$, which is to say it is the limit of $\varphi(1@\gamma,1@0)$ for $\gamma < \beta$. And $\Phi[n] = \varphi^*(1_{\beta[n]})$, and so $\Phi$ is the limit of $\varphi^*(1_\gamma)$ for $\gamma < \beta$. It is obvious that these limits are equal, due to the inductive hypothesis yielding $\varphi^*(1_\gamma) < \varphi(1@\gamma, 1@0) < \varphi^*(1_{\gamma+1})$. This method of showing two sets have the same suprema, by showing their elements are intertwined, is a useful method and a reader with a keen eye will notice we used it previously, in the proof of Theorem 3.10.
    \item If $\Phi$ is of class $C$, $\Phi' = \varphi(\alpha@0, \cdots)$ and $\Phi = \varphi^*(\alpha_0, \cdots)$
    \begin{enumerate}
        \item If $\alpha$ is a limit, then $\Phi'$ must be the limit of $\varphi(\gamma@0, \cdots)$ for $\gamma < \alpha$, as derivatives of normal functions are always normal. $\Phi$ is also the limit of $\varphi^*(\gamma_0, \cdots)$ for $\gamma < \alpha$. These are exactly equal.
        \item If $\alpha$ is a successor, by the above lemma (and the fundamental sequences) we have $\Phi$ being the next fixed point of some normal function $f$ after $\varphi^*(\dot{\alpha}_0, \cdots)$. Then, we can use an argument similar to cases 2, 3, 5, 6, and 7.
    \end{enumerate} 
    \item If $\Phi$ is of class $D$, then $\Phi' = \varphi(\cdots,\alpha@\beta)$, where $\alpha$ and $\beta$ are successors. In this case, $f((\cdots,\alpha@\beta)) = 2$, and $g((\cdots,\alpha@\beta),\xi) = (\cdots,\dot{\alpha}@\beta,\xi@\dot{\beta})$. Therefore, $\Phi'$ is the least fixed point of $\xi \mapsto \varphi(\cdots,\dot{\alpha}@\beta,\xi@\dot{\beta})$. And $\Phi = \varphi^*(\alpha_\beta, \cdots)$, where both $\alpha_\beta$ and $\beta$ are successors. The fundamental sequences imply that it is the least fixed point of $\xi \mapsto \varphi^*(\xi_{\dot{\beta}},\dot{\alpha}_\beta, \cdots)$. These are exactly equal.
    \item If $\Phi$ is of class $E$, then $\Phi' = \varphi(\cdots,\alpha@\beta)$, where $\alpha$ is a limit. Then, $f((\cdots,\alpha@\beta)) = \alpha$, and $g((\cdots,\alpha@\beta),\xi) = (\cdots,\xi@\beta)$.Then, $\Phi'$ is the least fixed point of $\xi\mapsto \varphi(\cdots,\gamma@\beta,\xi@0)$ for $\gamma<\alpha$, or the limit of $\varphi(\cdots,\gamma@\beta, 1@1)$ for $\gamma<\alpha$.
    $\Phi$ is $\varphi^*(\alpha_\beta, \cdots)$ for $\alpha_\beta$ a limit, and by the fundamental sequences it is the limit of $\varphi^*(\xi_\beta, \cdots)$ for $\xi < \alpha$. It is obvious that these limits are equal. 
    \item If $\Phi$ is of class $F$, then $\Phi = \varphi(\cdots,\alpha@\beta)$, where $\alpha$ is a successor and $\beta$ is a limit. Then, $f((\cdots,\alpha@\beta)) = \beta$, and $g((\cdots,\alpha@\beta),\xi) = (\cdots,\dot{\alpha}@\beta,1@\xi)$. Thus, $\Phi'$ is the least fixed point of $\xi\mapsto \varphi(\cdots,\dot{\alpha}@\beta,1@\gamma,\xi@0)$ for all $0<\gamma<\beta$\footnote{The restriction here is because otherwise the input of $\varphi$ is not in $A$.}, or the limit of $\varphi(\cdots,\dot{\alpha}@\beta,1@\gamma,1@1)$\footnotemark[5] for all $1<\gamma<\beta$.  $\Phi$ is $\varphi^*(\alpha_\beta, \cdots)$ for $\alpha_\beta$ a successor and $\beta$ a limit, and is the limit of $\varphi^*(1_\gamma,\dot{\alpha}_\beta)$ for $\gamma<\beta$. These limits are also equal. \qedhere
\end{enumerate}
\end{proof}

\begin{conjecture}
Let $\mathcal{C} = C_0(\mathcal{J}) \cap \Omega$. Also, let $+_{\mathrm{NNF}}$ (``notational normal form'') be the restriction of $+$ to $(s,t) \in \mathcal{OT}$ where $s + t \in \mathcal{OT}$ too, and $+_{\mathrm{ONF}}$ (``ordinal normal form'') be the restriction of $+$ to $(\alpha, \beta) \in \mathcal{C}$ where $\alpha + \beta$ is in Cantor normal form. Then there is an isomorphism

\begin{equation}
o: (\mathcal{OT}, \leq, +_{\mathrm{NNF}}) \longrightarrow (\mathcal{C}, \leq, +_{\mathrm{ONF}})
\end{equation}

\end{conjecture}

\begin{proof}[Proof Sketch]
Simultaneously define $q: \mathcal{A} \to A$ and $o: \mathcal{T} \to \mathrm{Ord}$ by:

\begin{enumerate}
    \item $o(\ulcorner 0 \urcorner) = 0$.
    \item $o(\ulcorner a+b \urcorner) = o(a)+o(b)$.
    \item $o(\ulcorner \varphi a \urcorner) = \varphi q(a)$.
    \item $q(\epsilon) = \emptyset$
    \item $q((\alpha_1 @ a_1, \alpha_2 @ a_2, \cdots, \alpha_n @ a_n)) = \{(o(\alpha_m), q(a_m)): 1 \leq m \leq n\}$.
\end{enumerate}

$o$ is obviously total and preserves addition, so we verify that:

\begin{enumerate}
    \item For all $s \in \mathcal{OT}$, $o(s) \in \mathcal{C}$.
    \item If $s, t \in \mathcal{OT}$ and $o(s) = o(t)$, then $s = t$.
    \item If $\alpha \in \mathcal{C}$, there is $s \in \mathcal{OT}$ so that $o(s) = \alpha$.
    \item For all $s, t \in \mathcal{OT}$, $s \leq t$ iff $o(s) \leq o(t)$.
\end{enumerate}

The first can likely be shown using an ``inverse'' algorithm for converting dimensional Veblen back into Buchholz's function. Such an algorithm has not been found yet, due to issues such as $\omega^{\varepsilon_0+1}$ corresponding to $\psi_0(\Omega+1)$ rather than $\psi_0(\varepsilon_0+1)$ -- however, it is still likely possible to find such an algorithm.

We define a different rank function, $\varrho: \mathcal{T} \to \mathbb{N}$, by:

\begin{enumerate}
    \item $\varrho(\ulcorner 0 \urcorner) = 0$.
    \item $\varrho(\ulcorner a+b \urcorner) = \max(\varrho(a), \varrho(b))+1$.
    \item $\varrho(\ulcorner \varphi (\alpha_1 @ a_1, \cdots, \alpha_n @ a_n) \urcorner) = \sup(\{\varrho(\alpha_i): 1 \leq i \leq n\} \cup \varrho''\mathbf{S}(a_1) \cup  \cdots \varrho''\mathbf{S}(a_n))+1$.
\end{enumerate}

In particular, $\varrho(\varphi \epsilon) = 1$. One then can prove injectivity by induction on $\max(\varrho(s), \varrho(t))$. If $\varrho(s) = \varrho(t) = 0$, then we have $o(s) = o(t) = 0$, and $s = t = 0$ too. In the induction step, one knows that that either $s$ or $t$ is nonzero, so it can be broken up into eight cases:

\begin{enumerate}
    \item When $s = 0$ and $t$ is nonprincipal.
    \item When $s = 0$ and $t$ is principal.
    \item When $s$ is nonprincipal and $t = 0$.
    \item When both $s$ and $t$ are nonprincipal.
    \item When $s$ is nonprincipal and $t$ is principal.
    \item When $s$ is principal and $t = 0$.
    \item When $s$ is principal and $t$ is nonprincipal.
    \item When both $s$ and $t$ are principal.
\end{enumerate}

All cases except the last two can be easily resolved, with the fact that $a+b$ being standard implies $o(a)+o(b)$ is in Cantor normal form, and $s$ being principal implies $o(s)$ being additively principal. For when $s$ and $t$ are both principal, one would need to analyse all the fixed point properties that they have, and show that standardness implies neither $s$ nor $t$ could be written in a simpler form -- e.g. neither can contain stuff like $\varphi(1,\varphi(2,0))$.

Proving that the comparison algorithm precisely coincides with real comparison would follow from a similarly detailed case classification. Surjectivity follows from the fact that $q$ is surjective and the algorithm converting Buchholz's function into dimensional Veblen.
\end{proof}

It actually turns out that $\mathcal{C} = \eta_0$ (a corollary of $C_\nu(\alpha) \cap \Omega_{\nu+1} = \psi_\nu(\alpha)$, proven in \cite{buchholz}). Therefore, if this conjecture is true then we have the following two corollaries:

\begin{corollary}
Assume $R$ is a recursive well-order on $\mathbb{N}$ so that transfinite induction along $R$ is provable in $\mathrm{ID}_1$, Peano arithmetic augmented by an axiom schema of non-iterated inductive definitions, or $\mathrm{KP}$, Kripke-Platek set theory with the axiom of infinity. Then $(\mathbb{N}, R)$ is isomorphic to an initial segment of $(\mathcal{OT}, \leq)$.
\end{corollary}

\begin{proof}
The ordinal analysis of $\mathrm{ID}_1$ was given in \cite{buchholz3}, and there have been many ordinal analyses of $\mathrm{KP}$ -- one can be found in \cite{jager}. Both theories have the same proof-theoretic ordinal $\eta_0$, which implies that every recursive well-order on $\mathbb{N}$ with provable transfinite induction is isomorphic to $(\eta_0, \in)$ or an initial segment. And this is itself isomorphic to $(\mathcal{C}, \in)$.
\end{proof}

\begin{corollary}
There is a nonstandard model $\mathcal{M} = (M, E)$ so that $\mathcal{M} \models \mathrm{KP} + \ulcorner (\mathcal{OT}, \leq)$ is ill-founded $\urcorner$.
\end{corollary}

\begin{proof}
By consistency of $\mathrm{KP}$ and the completeness theorem.
\end{proof}

This precise analysis of the strength of $\varphi$ also yields a system of fundamental sequences below $\eta_0$, which could be used to build a Hardy or fast-growing hierarchy (cf. \cite{buchholz2}) Namely, since $o$ is a bijection, if one lets $L = \varrho \circ o^{-1}$ then

\begin{equation}
L: \eta_0 \to \mathbb{N}
\end{equation}

is total and we claim that $L$ is a norm -- this notion was introduced in \cite{buchholz4}:

\begin{proposition}
$L$ is a norm.
\end{proposition}

\begin{proof}
We first show $L(0) = 0$: $L(0) = \varrho(o^{-1}(0)) = \varrho(\ulcorner 0 \urcorner) = 0$. Next, we need to show that, for $\alpha < \eta_0$, we have $L(\alpha+1) \leq L(\alpha)+1$. Let $\alpha = o(a)$, where $a \in \mathcal{OT}$. If $a \in \mathcal{PT}$ then $\ulcorner a + \varphi \epsilon \urcorner \in \mathcal{OT}$, $o(\ulcorner a + \varphi \epsilon \urcorner) = \alpha+1$ and $\varrho(\ulcorner a + \varphi \epsilon \urcorner) = L(\alpha)+1$, so $L(\alpha+1) = L(\alpha)+1$. If $a \notin \mathcal{OT}$, let $a = a_1 + a_2 + \cdots + a_n$ where all the $a_i$'s except $a_n$ are principal. Then $\ulcorner a_1 + a_2 + \cdots + a_n + \varphi \epsilon \urcorner \in \mathcal{OT}$ (where non-binary sums are intended to be right-associative), $o(\ulcorner a_1 + a_2 + \cdots + a_n + \varphi \epsilon \urcorner) = \alpha+1$ and $\varrho(\ulcorner a_1 + a_2 + \cdots + a_n + \varphi \epsilon \urcorner) = L(\alpha)+1$, so $L(\alpha+1) = L(\alpha)+1$. Lastly, we need to prove that, for each $n < \omega$ and $\alpha < \eta_0$ we have $|\{\beta < \alpha: L(\beta) < n\}| < \aleph_0$. This is straightforward to show, as the set of terms with rank at most any fixed natural number is always finite, since each element of the next stage is generated from the previous one in a very specific, finite amount of ways.
\end{proof}

Actually, $L$ commutes with Cantor normal form like so (remember, this is still assuming Conjecture 3.16): if $\alpha$ has CNF $\omega^{\alpha_0} + \omega^{\alpha_1} + \cdots + \omega^{\alpha_n}$ and $n > 0$ then $L(\alpha) = \max(L(\omega^{\alpha_0}), L(\omega^{\alpha_1}), \cdots, L(\omega^{\alpha_n})) + n$ and, if $\alpha \neq \omega^\alpha$ (i.e. $\alpha$ is not an epsilon number) then $L(\omega^\alpha) = L(\alpha) + 1$. In particular, if $\alpha < \varepsilon_0$ has CNF $\omega^{\alpha_0} + \omega^{\alpha_1} + \cdots + \omega^{\alpha_n}$ then

\begin{equation}
L(\alpha) = \max(L(\alpha_0), L(\alpha_1), \cdots, L(\alpha_n)) + n + 1
\end{equation}

As such, this induces a fundamental sequence system $\alpha[n] = \max\{\beta < \alpha: L(\beta) \leq L(\alpha) + n\}$ for $\alpha < \eta_0$. However, past the limit of the Wainer hierarchy (i.e. $\varepsilon_0$) this fundamental sequence system would likely have quite curious behaviour, not agreeing with the fundamental sequences below $E_0$ introduced earlier, and is just a specific instance of a wider framework that has already been investigated in detail.

As some closing remarks, we would like to note that, beyond the level of the Large Veblen Ordinal, Buchholz's $\psi$-functions are vastly superior in terms of utility and ease of definition, considering how long the definition of, and proofs regarding, the dimensional Veblen function are. However, we believe that the dimensional Veblen function still has merit, because it can be considered as taking the idea of iterated fixed points to its limit, and not stopping when the limit of 1-dimensional structures is reached. The dimensional Veblen function may also be substantially useful as a benchmark for analysing the strength of other ordinal representation systems.

On another note, it is possible to consider $\mathrm{ID}_1$ to consist of only first-order arithmetic, with increased induction, and $\mathrm{KP}$ to be precisely the predicative or recursive fragment of ZFC (e.g. $\alpha$-recursion theory studies generalized computability theory within models of $\mathrm{KP}$). Therefore, these theories could be considered the limit of predicative reasoning, and dimensional Veblen as a notation exhausting all the predicative ordinals. This makes sense, as ordinal-representation systems beyond $\eta_0$ (e.g. for admissible set theory, second-order arithmetic and iterated inductive definitions) all heavily make use of uncountable or nonrecursive ordinals. This may be considered a bit of a stretch because, as per \cite{weaver}, one may consider the true limit of predicative ordinals (that have been found so far) to be the small Veblen ordinal, which is significantly smaller than $\eta_0$. Also, our notion of predicative -- ordinals with associated recursive well-orders, transfinite induction along which is provable in first-order arithmetic with strengthened induction or set theory without powerset and recursive separation and predicative collection -- doesn't quite align with Weaver and Feferman's notion -- ordinals with associated recursive well-orders, well-orderedness which is provable in predicative fragments of second-order arithmetic. This is discussed in more detail in \cite{weaver}.

Also, considering some extensions, we may think of $(1@[1,0])$ as the ``least'' (under $\prec$) fixed point of $X \mapsto (1@X)$, and then think of $(1@[1,1])$ as the second, $(1@[2,0])$ as the fixed point of $X\mapsto (1@[1,X])$, and so on. Further extensions (once formalized) with more types of brackets may reach the proof-theoretical ordinal of $\Pi^1_1\mathsf{-CA_0}$. Having higher brackets as analogues of higher uncountables in Buchholz's $\psi$, it can easily be extended to use cardinals such as $\Omega_{\omega+1}$, $\Omega_\Omega$, and $\Omega_{\Omega_\Omega}$. This ``extended Buchholz's $\psi$'' may be extended to have 3 arguments or more, as the Veblen hierarchy $\varphi_\alpha$ can be extended to have 3 or more arguments. Also of some interest is the fact that $\psi$'s $\Omega$ roughly corresponds to the arguments and dimensions in $\varphi$, and this chain could be further extended, with a new function whose $\Omega$ corresponds to the arguments in $\psi$. It is possible that this chain could extend beyond $(^{++})$-stability -- which would yield a new paradigm of ordinal notations, as many in the literature are very complicated. However, such an extension would likely be too synthetic.

\printbibliography[heading=bibintoc,title={References}]
\end{document}